\documentclass[12pt]{article}
\usepackage{amssymb}
\usepackage{MnSymbol}
\usepackage{amsmath}
\usepackage{amsbsy}
\usepackage{latexsym}
\usepackage{graphicx}
\usepackage[all]{xy}

\newtheorem{cor}{Corollary}[section]
\newtheorem{defn}{Definition}[section]
\newtheorem{hyp}{Hypothesis}[section]
\newtheorem{lemma}{Lemma}[section]
\newtheorem{prop}{Proposition}[section]

\newtheorem{rem}{Remark}[section]
\newtheorem{thm}{Theorem}[section]

\newcounter{example}[section]
\newcommand{\addwords}[1]{\ \ \ \mbox{#1}}

\newcommand{\ann}{{\rm ann}}
\newcommand{\Aut}{{\rm Aut}}
\newcommand{\ba}{\ensuremath{{\bf a}}}

\newcommand{\PVmodbarG}{\ensuremath{
\bar{G}\backslash \bbP V}}

\newcommand{\bc}{\ensuremath{{\bf c}}}
\newcommand{\bbC}{\ensuremath{\mathbb{C}}}

\newcommand{\bbN}{\ensuremath{\mathbb{N}}}
\newcommand{\bbP}{\ensuremath{\mathbb{P}}}
\newcommand{\bbQ}{\ensuremath{\mathbb{Q}}}

\newcommand{\bbR}{\ensuremath{\mathbb{R}}}
\newcommand{\bbT}{\ensuremath{\mathbb{T}}}
\newcommand{\bbZ}{\ensuremath{\mathbb{Z}}}

\newcommand{\bD}{\ensuremath{{\bf D}}}
\newcommand{\be}{\ensuremath{{\bf e}}}
\newcommand{\beqar}{\begin{eqnarray}}	
\newcommand{\beq}{\begin{equation}}
\newcommand{\beql}[1]{\begin{equation}\label{#1}}

\newcommand{\bPf}{\noindent \textsc{Proof\ }}

\newcommand{\bv}{\ensuremath{{\bf v}}}
\newcommand{\bu}{\ensuremath{{\bf u}}}
\newcommand{\bw}{\ensuremath{{\bf w}}}
\newcommand{\bx}{\ensuremath{{\bf x}}}

\newcommand{\cA}{\ensuremath{\mathcal{A}}}
\newcommand{\cB}{\ensuremath{\mathcal{B}}}

\newcommand{\cD}{\ensuremath{\mathcal{D}}}
\newcommand{\cE}{\ensuremath{\mathcal{E}}}
\newcommand{\cF}{\ensuremath{\mathcal{F}}}
\newcommand{\cG}{\ensuremath{\mathcal{G}}}
\newcommand{\cf}{\emph{cf.}}
\newcommand{\cH}{\ensuremath{\mathcal{H}}}
\newcommand{\cI}{\ensuremath{\mathcal{I}}}

\newcommand{\cM}{\mathcal{M}}

\newcommand{\cO}{\ensuremath{\mathcal{O}}}

\newcommand{\cR}{\ensuremath{\mathcal{R}}}

\newcommand{\cT}{\ensuremath{\mathcal{T}}}
\newcommand{\cV}{\ensuremath{\mathcal{V}}}

\newcommand{\cY}{\ensuremath{\mathcal{Y}}}

\newcommand{\displaymapdef}[5]
{
	\begin{array}{rcrcl}
		#1 &:& #2 &\longrightarrow& #3 \\
		& &    &                    \\
		& & #4 &\longmapsto    & #5
	\end{array}
}
\newcommand{\displaynomapdef}[4]
{
	\begin{array}{ccc}
		#1 &\longrightarrow& #2 \\
		&                    \\
		#3 &\longmapsto    & #4
	\end{array}
}

\newcommand{\eeq}{\end{equation}}
\newcommand{\eeqar}{\end{eqnarray}}
\newcommand{\eg}{\emph{e.g.}}
\newcommand{\eitherortw}[4] 
{\left\{                  
	\begin{array}{ll}         
		#1& \mbox{#2}\\
		#3& \mbox{#4}
	\end{array}
	\right.}
\newcommand{\eitherorth}[6] 
{\left\{                  
	\begin{array}{ll}         
		#1& \mbox{#2}\\
		#3& \mbox{#4}\\
		#5& \mbox{#6}
	\end{array}
	\right.}
\newcommand{\End}{{\rm End}}
\newcommand{\ePf}{~$\Box$\vertsp\par}
\newcommand{\eps}{\varepsilon}

\newcommand{\fa}{\ensuremath{\mathfrak{a}}}

\newcommand{\fD}{\ensuremath{\mathfrak{D}}}
\newcommand{\ff}{\ensuremath{\mathfrak{f}}}
\newcommand{\fg}{\ensuremath{\mathfrak{g}}}

\newcommand{\fl}{\ensuremath{\mathfrak{l}}}

\newcommand{\fz}{\ensuremath{\mathfrak{z}}}

\newcommand{\GL}{{\rm GL}}
\newcommand{\half}{\frac{1}{2}}
\newcommand{\Hom}{{\rm Hom}}
\newcommand{\id}{{\rm id}}
\newcommand{\ie}{{\em i.e.}}
\newcommand{\im}{{\rm im}}
\newcommand{\Inn}{{\rm Inn}}
\newcommand{\inv}{^{-1}}
\newcommand{\ndiv}{\nmid}
\newcommand{\op}{^{\rm op}}

\newcommand{\Out}{{\rm Out}}

\renewcommand{\rem}{\refstepcounter{rem}\noindent{\sc Remark \therem}}

\newcommand{\resp}{{\em resp.}}

\newcommand{\ses}[5]{#1\rightarrow #2\longrightarrow #3\longrightarrow 
#4\rightarrow #5}

\newcommand{\SL}{{\rm SL}}
\newcommand{\smallcolvec}[2] 
{
\left( \begin{smallmatrix}
	#1\\
	#2
\end{smallmatrix} \right)
}
\newcommand{\smalllegsym}[2] 
{
\left( \begin{smallmatrix}
	#1\\
	-\\
	#2
\end{smallmatrix} \right)
}
\newcommand{\smallttmat}[4] 
{
	\left( \begin{smallmatrix}
		#1&#2\\ #3&#4
	\end{smallmatrix} \right)
}

\newcommand{\Sp}{{\rm Sp}}

\renewcommand{\t}[1]{\tilde{#1}}

\renewcommand{\therem}{\thesection.\arabic{rem}}

\newcommand{\Tr}{{\rm Tr}}

\newcommand{\U}{{\rm U}}
\newcommand{\ut}{\underline{t}}
\newcommand{\vertsp}{\vspace{1ex}}

\newcommand{\Zover}[1]{\ensuremath{\bbZ/#1}}
\newcommand{\Zunder}[1]{#1 \inv\bbZ/\bbZ}
\newcommand{\AdsB}{\ensuremath{A\oplus B}}
\newcommand{\AIdsB}{\ensuremath{A_{I,\ff}\oplus B_{I,\ff}}}
\newcommand{\barg}{\ensuremath{\bar{g}}}
\newcommand{\barG}{\ensuremath{\bar{G}}}
\newcommand{\barphi}[1]{\bar{#1}_{\bar \phi}}
\newcommand{\cAng}{{\frak O}}

\newcommand{\Iproddots}{\ensuremath{\langle \cdot, \cdot\rangle}}
\newcommand{\Hei}[4]{\cH({#1},{#2},{#3},{#4})}
\newcommand{\HeiStan}{\Hei{A}{B}{C}{\lambda}}

\newcommand{\Of}{{\cO/\ff}}
\newcommand{\Poly}{{\rm Poly}}

\newcommand{\rep}{representation\ }
\newcommand{\SpABd}[1]{\ensuremath{{\rm Sp}_{#1}(\AdsB;\delta)}}

\begin{document}
\title{Towards a Theory of SIC-like Phenomena:\\ 
Regular Bouquets and Generalised Heisenberg Groups}
\author{David Solomon}
\date{\today}
\maketitle
\begin{abstract} 
	\noindent We lay the foundations for a broad algebraic theory encompassing SICs in the hope of elucidating their heuristic connections with Stark units. What emerges is a greatly generalised set-up with added structure and potential for applications in other areas. Let $A$ and $B$ be finite modules for a commutative ring $R$, $C$ a finite abelian group and $\lambda: A\times B\rightarrow C$ an $R$-balanced bilinear pairing. The main constructs are the generalised Heisenberg group $\cH=\HeiStan$ attached to these data (an abstract central extension of $A\oplus B$ by $C$) which plays the role of the Weyl-Heisenberg group in SIC theory,  together with its canonical, unitary Schr\"odinger representations.
	The SIC itself is replaced by an $\cH$-orbit of complex lines  in the representation space, termed a \textit{bouquet}. The overlaps of the SIC are interpreted as a map from $\cH/Z(\cH)$ into $\bbC$ whose absolute values are the Hermitian `angles' between the lines in the bouquet. We also introduce a regularity condition on bouquets in terms of the angle-map, intended to weaken the equiangularity condition of SICs. At the same time, it allows the incorporation of the $R$-structure \textit{via} the abstract automorphism group of $\cH$ which in turn generalises the Clifford group of SIC theory \textit{via} its associated Weil representation. As well as several subsidiary definitions and `structural' results, we prove a new `clinometric relation' for the angle-map, determine the structure of the automorphism group and introduce a large class of examples of arithmetic origin, derived from the trace-pairing on quotients of fractional ideals in an arbitrary number field, which we investigate in greater detail. 
\end{abstract}
\section{Introduction}
\subsection{Background}
This paper is inspired by the remarkable results and observations made by many different workers in the theory of SICs (or SIC-POVMS) over the past 25 years, and particularly by the number-theoretic connections that they have more recently uncovered. 

Depending on applications and context, the idea of a `SIC' is interpreted in different but equivalent ways in the literature. Two common interpretations are as maximal sets of equiangular lines in $\bbC^s$ (with $s>3$) and as certain sets of Hermitian $s\times s$ matrices. We briefly explain the first viewpoint which is the one adopted here. (For the second, see for instance~\cite{Kopp IMRN}, and for the correspondence between the two viewpoints in an equivariant context, see \S\ref{subsec:Projectors}.) 

For us, a SIC is simply a set of $s^2$ complex lines through the origin in $\bbC^s$ which is \textit{equiangular}   
in the sense that, the `angle' between two unit vectors generating any pair of distinct lines (\ie\ the absolute value of their standard Hermitian inner product) is independent of the pair chosen. It was shown in~\cite{Delsarte et al} that no more than $s^2$ equiangular lines can exist in $\bbC^s$, but it it is still unknown whether this bound is attained (\ie\ whether SICs exist) for every $s>3$. One part of \textit{Zauner's Conjecture} affirms that this is so, see~\cite{Zauner}.  

Serious interest in SICS \textit{per se} began at the tail end of the last century, in  the areas of Quantum Designs and, later, Finite Tight Frames (see Remark~\ref{rem:tight frames}). They find applications in signal processing, orthogonal polynomials, quantum information theory~\textit{etc}. (The last of these is responsible for much of the SIC-nomenclature used here, starting with the term `SIC-POVM': an abbreviation of \textit{Symmetric, Informationally Complete Positive Operator-Valued Measure}.) There is by now a large literature concerning the theory and applications of SICs, comprising both proven results and extensive computations that give rise to heuristic observations, known as \textit{`SIC phenomenology'}. For a sampling of this work, see the papers 
\cite{ABGHM},
\cite{Appleby Flammia McC Yard}, 
\cite{DMA: SICPOVMs and the ECG}, 
\cite{Gal Auts Appleby et al}, 
\cite{Bengtsson Grass McC}, 
\cite{Horodecki et al}, 
\cite{Kopp IMRN}, 
and \cite{Zauner}. 
Waldron's book~\cite{Waldron's book} on tight frames gives a thorough introduction to the theory of SICs in Chapter~14.

A first and fundamental example of SIC phenomenology is the observation that every known SIC in $\bbC^s$ (except the \textit{Hoggar lines} for $s=8$) is a unitary transformation of a  \textit{`Heisenberg (equivariant) SIC'} namely the orbit of a single line under the discrete Weyl-Heisenberg group which we denote here ${\rm WH}(s)$. The latter is a certain finite subgroup of the group $\U(\bbC^s)$ of all unitary transformations of $\bbC^s$, which is nilpotent of class $2$ and of order $s^3$ if $s$ is odd and $2s^3$ if $s$ is even. Its image in ${\rm PU}(\bbC^s):={\rm U}(\bbC^s)/\bbC^\times$ has order $s^2$ in both cases. For details, we refer to \S\ref{subsec:the Base Case} where this appears as the simplest possible special case -- we call it the `Base Case' -- of a vastly generalised set-up introduced in \S\S\ref{reps of SV-type}--\ref{sec:Aut Gps}  . (More on the latter below.)

Increasingly though,  SICs and `SIC-like phenomena' are also catching the attention of number theorists. This is due to unexpected connections with units in the rings of integers of certain algebraic number fields. More precisely, in every one of the very many Heisenberg SICs computed so far, it has also been observed that the values of Hermitian inner products of the unit vectors (\textit{before} taking absolute values) are of form $(s+1)^{-1/2}$ multiplied 
by algebraic complex numbers of absolute value $1$, namely certain rational powers of Galois conjugates of \textit{Stark units} (see below) in the ray-class field $k_s((s)\infty_1\infty_2)$ over the real quadratic field $k_s:=\bbQ(\sqrt{(s-3)(s+1)})$. The reader may refer, for instance, to~\cite{ABGHM} for an account of this heuristic connection and to~\cite{Kopp IMRN} for a precise conjecture in a special case. 
Very recently, the first, non-heuristic `explanation' of this connection has appeared in~\cite{AFK constructive}. Complex-analytic in nature, it makes use of the \textit{Shintani-Fadeev cocycle} and relies on two conjectural results, including the (abelian, first order) Stark Conjecture itself over the fields $k_s$. 

Now, on the one hand, the Stark Conjecture affirms the existence of certain, essentially unique algebraic units (the Stark units) in abelian extension fields of $k_s$, the logarithms of whose real absolute values are specified in terms of the derivatives at $z=0$ of certain Artin $L$-functions proper to $k_s$. (See~\cite{Stark}, \cite{Tate} for more details.) 
On the other hand, if they do exist, Stark himself showed that in many cases the Stark units generate the abelian extensions in question, over $k_s$. This gives rise to the hope that a better understanding of SICs and related objects could lead to advances in the study of Stark's conjecture and/or that of Hilbert's 12th problem (the analytic construction of class-fields) over such real quadratic fields or even, perhaps, more generally.
\subsection{Overview of this Paper}
Our principal aim here is to lay the foundations of a broad \textit{algebraic} theory around SICs that, by adding several new layers of structure and generality, may help to fulfil the above hope. 
In fact, the algebraic structure we uncover here seems sufficiently general and extensive to be of interest in itself and, maybe, even to lead to new, SIC-like phenomena in other areas of mathematics. 

This theory can be described in the vaguest and most general terms as  \textit{`a new  machine for attaching systems of numbers to a certain class of algebraic structures'.} More precisely, the class consists of tuples  $(A,B,C,\lambda)$ (plus certain other data) where $A$ and $B$  are finite (later profinite) modules for a commutative ring $R$, $C$ is a finite abelian group and $\lambda:A\times B\rightarrow C$ is an $R$-balanced, $\bbZ$-bilinear form. Such tuples can be seen as the objects of a category whose morphisms are triples $(t_A,t_B,t_C)$ of homomorphisms satisfying certain conditions (see Proposition~\ref{prop: criterion for diagonal homs} and preceding). From this viewpoint, we are concerned here  with certain functors from this category \eg\ a `Schr\"odinger functor' to the category of unitary representations of finite groups. Rather than pursue the categorical viewpoint here, however, we now describe informally several interconnected `cogs' of the machine which we shall put together piece-by-piece in Sections~\ref{reps of SV-type}--\ref{sec:Aut Gps}.

 The drive-wheel of this machine is the  Generalised Heisenberg Group (or GHG) $\cH=\HeiStan$ attached to the tuple $(A,B,C,\lambda)$. This is a certain abstract, finite group of nilpotency class $\leq 2$. 
The $R$-structure on $A$ and $B$ is a novel feature introduced here to enrich the theory. It does not affect the group structure of $\cH$ itself but it does  provide an extra mechanism for encoding  algebraic information, as we shall later see. Without it (or by taking $R$ to be $\bbZ$, the default option) our definition of GHGs coincides with~\cite[Def.\ 4.1]{Szabo} and generalises the matrix group ${\rm WH}(s)$. (For other versions of finite Heisenberg groups, see for example~\cite{Tata III}, \cite{Howe} and~\cite{Wi: Heisenberg Gp}.) In fact, Szab\'o shows  that \textit{every} class $\leq 2$ nilpotent group is a subgroup of some $\HeiStan,$ even with the added condition that $\lambda$ be non-degenerate. This condition, together with cyclicity of $C$,  is in fact our normal working Hypothesis~ND-C for GHGs. It ensures, for instance, that the quotient of $\cH$ by its centre $Z$ is (canonically) isomorphic to $A\oplus B$ and hence inherits an $R$-module structure.

Connected directly to the drive-wheel are, on the one hand, the subgroup $\Aut_R^0(\cH)$ of $\Aut(\cH)$ (automorphisms which fix $Z$ elementwise and are $R$-linear modulo $Z$) and on the other, a pair of canonical `left and right \textit{Schr\"odinger representations} of $\cH$ which we denote $\sigma_p$ and $\tau_p$. They are siomorphic so we restrict attention largely to $\sigma_p$. Its representation space is denoted $\cM(A)$, has dimension $s=|A|$. It is unitary and  (under Hypothesis~ND-C plus a mild condition) irreducible and of \textit{`Stone-Von Neumann Type'} (`SV' for short) in the terminology of~\cite{Howe}. In the Base Case with $s$ odd, ${\rm WH}(s)$ is essentially the image in the unitary group ${\rm U}(\bbC^s)$ of $\cH(\Zover{s},\Zover{s},\Zover{s}, \times)$ under the faithful, left or right Schr\"odinger representation (see Subsection~\ref{subsec:the Base Case}). 

An important part of the mechanism connecting these two cogwheels is the (projective) \textit{Weil representation} of $\Aut^0(\cH)=\Aut^0_\bbZ(\cH)$ associated to $\sigma_p$ when the latter is irreducible. In the Base Case, for example, the Weil representation (essentially) embeds $\Aut^0(\cH)$ in the projective unitary group ${\rm PU}(\bbC^s)$. The image is known in the Physics literature as the \textit{`Clifford group'}. It has been observed heuristically to `carry a piece of the Galois action' on the associated Stark units in this case (see~\cite[Ch.\ 14]{Waldron's book} or~\cite{ABGHM}, for example).

The Weil representation also allows $\Aut^0(\cH)$ to act (via its quotient $\Out^0(\cH)$) on another important component of the machine, namely the set of `$\cH$-bouquets'. Each such bouquet $\cY$ is defined to be an $\cH/Z$-orbit of (complex) lines in $\cM(A)$. It comes equipped with an \textit{overlap-map} from $\cH/Z$ into $\bbC$ whose absolute values constitute the so-called \textit{angle-map} because they represent the complex angles between unit generators of the  lines. We show that the squares of these absolute values must satisfy a certain `clinometric relation' (Theorem~\ref{thm:clinomteric thm for nilp gps}). Moreover, if $\cY$ is free, its cardinality is $s^2$ so it will be an $\cH$-equivariant SIC in $\cM(A)$ iff the angle-map is constant on the  non-identity elements of $\cH/Z$.

In \S\ref{sec:Overlap- and Angle-Maps, and Bouquets} we introduce a novel, and perhaps more natural condition on bouquets which weakens equiangularity. This is  \textit{$\cA$-regularity} where $\cA$ is any subgroup of $\Out^0(\cH)$ (for instance the group $\Out^0_R(\cH)$ which is the image of $\Aut^0_R(\cH)$). Conceptually, it means that the unit vectors of $\cY$ `span an $\cA$-regular complex polytope' in a sense that we shall not make completely formal here. 
We mention in passing another condition on bouquets, namely \textit{invariance}, introduced in~\S\ref{sec: nilpotent groups} to generalise yet another aspect of SIC phenomenology: the symmetry of order $3$ discovered by Zauner.    

As for the input tuples $(A,B,C,\lambda)$ to be fed into the machine, we give a large class of arithmetic examples in \S\ref{sec:Examples} derived from the trace-map on fractional ideals modulo an integral ideal $\ff$, both w.r.t.\ the maximal order $\cO_k$ of a number field $k$. The ring $R$ above is taken to be $\cO_k/\ff$. (In principle, this construction could be generalised to non-maximal orders of $k$.) We refer to the associated groups $\cH=\cH(A,B,C,\lambda)$ as GHGs \textit{of arithmetic type}. It seems to the author that $\Out^0_R(\cH)$-regular bouquets for such $\cH$ form a natural context for a generalised theory of Heisenberg SICs in which to investigate their phenomenology. In fact, Heisenberg SICs are simply equiangular bouquets in the case $k=\bbQ$. This case is discussed in detail in \S\ref{subsec:the Base Case}.

A feature of the general arithmetic examples is that the modules $A$ and $B$ are \textit{free of rank one} over $R$, although usually without canonical generators. This leads to a proof that in such cases (and when $|\cH|$ is odd) $\Aut^0_{\cO_k/\ff}(\cH)$ is  isomorphic to 
a certain semidirect product of $\SL_2(\cO_k/\ff)$  (Theorem~\ref{thm: 3x2 diagram of aut gps} plus Theorem~\ref{thm: Sp iso to SL in general arithmetic case}). This generalises a result of~\cite{DMA: SICPOVMs and the ECG} in the Base Case. The proof of Theorem~\ref{thm: 3x2 diagram of aut gps}  (concerning $\Aut^0_{R}(\cH)$ for \textit{any} GHG $\cH$ with $R$-structure,  satisfying Hypothesis~ND-C and of odd order) relies in turn on a natural generalisation of the `classical'  \textit{displacement operators} of SIC-theory that we introduce in \S\ref{sec:Aut Gps}.   

\subsection{Perspectives}
Some further lines of research for GHGs of arithmetic type, including their $p$-adic limits and representations, are outlined in~\S\ref{subsec:perspectives}. It would also be interesting to have a classification of their regular bouquets in simple situations. The discussion in~\S~\ref{subsec:the Base Case} suggests that, even in the Base Case, further conditions on regular bouquets are required before they can be expected to produce algebraic generalisations of the specific Stark units heuristically associated to SICs. Another intriguing perspective is to generalise the `$r$-SICs' developed in~\cite{AFK constructive} for the Base Case, by investigating higher-dimensional subspaces of $\cM(A)$ under the action of GHGs via their Schr\"odinger representations. 

A limiting factor in much of this research is our ability 
to perform calculations for numerical experimentation. 
Currently this has been done only for SICs in the Base Case, by Scott, Grassl, Flammia and others, and already there it has proved to be very  demanding of computing power and time. Adapting their methods to investigate regular bouquets, for example, is probably fairly straightforward. However the computations required to treat many cases of significant interest are likely to be extremely intensive without substantial advances in methodology. 

Finally, and more distantly, the general machine developed here could be applied largely without change to  other natural classes of input tuples $(A,B,C;\lambda)$. For example, it would be interesting to apply it to tuples of algebro-geometric origins, especially in situations where $A$ and $B$ have more complicated $R$-structure.  
\subsection{Acknowledgements}
I am very pleased to be able to thank Marcus Appleby for many useful and inspiring conversations about SICs and for his expert guidance to their voluminous literature. I am grateful  also to I.H.E.S.\ where this work was started and partially written in March-May 2024. 
\subsection{Notations}  
If $M$ is a module for a commutative ring $R$, we shall write 
$\End_R(M)$ for $\Hom_R(M,M)$ considered as an $R$-algebra, with composition as multiplication, and $\Aut_R(M)$ for its unit-group. If $R$ is $\bbZ$, it will often be dropped from the notation. Thus $\Aut(A)$ simply denotes the automorphism group of an abelian group $A$ and we extend this notation to any group $G$. In this case, $\Inn(G)$ will denote the normal subgroup of $\Aut(G)$ consisting of all \textit{inner} automorhisms, and $\Out(G)$ will denote the quotient group $\Aut(G)/\Inn(G)$ of \textit{`outer automorphisms'}. If $X$ is an abelian group (\resp\ a ring) and $d$ is a positive integer, we shall write simply $X/d$ for the quotient group (\resp\ ring) $X/dX$. For any complex number $z$, we shall write $z^\ast$ for the complex conjugate and we write $\bbT$ for the circle group 
$\{z\,:\,z^\ast z=1\}<\bbC^\times $. (NB: its elements are usually referred to as \textit{`phases'} in the physics literature on SICs.) Its unique subgroup of order $f\in\bbN$ will be denoted $\mu_f$. 
A representation of a group $G$ on a non-zero vector space $V$ will be thought of either as a homomorphism $\rho:G\rightarrow\GL(V)$ or, equivalently, as a linear, left group action `$\cdot_\rho$' of $G$ on $V$. Thus $g\cdot_\rho v$ (or just $g\cdot v$) will mean the same thing as $\rho(g)(v)$, for any $g\in G$ and $v\in V$. We shall move between the two forms for purposes of  clarification or notational simplicity. Quotients by normal subgroups will be written $G/N$ and quotient sets (without group structure) by a group action will generally be written $H\backslash S$. 

\section{SV Representations and their Weil Representations}\label{reps of SV-type}
We begin by characterising the class of representations of finite groups, called {\em representations of Stone-Von Neumann type} in~\cite[\S3]{Howe} ({\em q.v.} for the reason for the name).
Let $G$ be finite group, $V$ a complex vector space of dimension $s$ and $\rho:G\rightarrow \GL(V)$ a representation of $G$ with character $\chi_\rho$. For any $g\in G$, the  eigenvalues of $\rho(g)$ are roots of unity, so $|\chi_\rho(g)|\leq s$. 
We recall that if $\rho$ is irreducible then the centre $Z=Z(G)$ acts by scalars through $\rho$ (by Schur's Lemma) so there is a (unique) homomorphism $\psi=\psi_\rho : Z\rightarrow \bbC^\times$, called  the {\em central character of $\rho$} such that $\rho(z)=\psi(z)\id_V$ for all $z\in Z$. 
In particular, $|\chi_\rho(z)|=s$ for all  $z\in Z$.  
\begin{thm}\label{thm:equivalent condits for SV} Suppose $\rho:G\rightarrow \GL(V)$ (not {\em a priori} irreducible) $\chi_\rho$, $s$ and $Z$ are as above and let $t$ denote $|G/Z|$. The following are  equivalent. 
	\begin{enumerate}
		\item\label{part1:equivalent condits for SV}  $\chi_\rho(g)=0$ for all $g\nin Z$, and $t=s^2$
		\item\label{part2:equivalent condits for SV}   $\rho$ is irreducible and $t=s^2$ 
		\item\label{partextra:equivalent condits for SV} $\rho$ is irreducible and $\chi_\rho(g)=0$ for all $g\nin Z$
		\item\label{part3:equivalent condits for SV} $\rho$ is irreducible and $\rho(g_1), \rho(g_2),\ldots ,\rho(g_t)$ are linearly independent in $\End_\bbC(V)$, for some set of representatives $g_1,g_2,\ldots,g_t$ of $Z$ in $G$
		\item\label{part4:equivalent condits for SV}
		$\rho$ is irreducible and $\rho(g_1), \rho(g_2),\ldots ,\rho(g_t)$ form a basis of $\End_\bbC(V)$, for any set of representatives $g_1,g_2,\ldots,g_t$ of $Z$ in $G$
	\end{enumerate}	
\end{thm}
\bPf 
Suppose $\rho\cong\bigoplus_{i=1}^n l_i\rho_i$ where the $\rho_i$ are irreducible and $l_i\in \bbN\,\,\forall i$. By character theory: 
$$
\sum_{i=1}^n l_i^2 = 
\frac{1}{|G|}\sum_{g\in G} |\chi_\rho(g)|^2=
\frac{1}{|G|}\sum_{g\in Z} |\chi_\rho(g)|^2+
\frac{1}{|G|}\sum_{g\nin Z} |\chi_\rho(g)|^2
$$
If~\ref{part1:equivalent condits for SV} holds, the R.H.S.\ equals $s^2 |Z|/|G|+0=1$, hence also the L.H.S., so $\rho$ is irreducible and~\ref{part2:equivalent condits for SV} and~\ref{partextra:equivalent condits for SV} must hold. 
If $\rho$ is irreducible the above equations give
\begin{equation}\label{eq:1=s2/r+...}
	1=\frac{s^2}{t}+
	\frac{1}{|G|}\sum_{g\nin Z} |\chi_\rho(g)|^2
\end{equation}
so~\ref{partextra:equivalent condits for SV} implies~\ref{part2:equivalent condits for SV} and~\ref{part2:equivalent condits for SV} implies~\ref{part1:equivalent condits for SV}. 
Recall that irreducibility also implies that 
$\{\rho(g)\}_{g\in G}$ spans $\End_\bbC(V)$ over $\bbC$, hence also  $\{\rho(g_1), \rho(g_2),\ldots ,\rho(g_r)\}$ for any set or representatives of $Z$ in $G$, since $Z$ acts by (non-zero) scalars. The equivalence of~\ref{part2:equivalent condits for SV}, \ref{part3:equivalent condits for SV} and \ref{part4:equivalent condits for SV} follows, since $\dim(\End_\bbC(V))=s^2$.
\ePf
\noindent Any representation $\rho:G\rightarrow \GL(V)$ satisfying the equivalent conditions of Theorem~\ref{thm:equivalent condits for SV} will be called an {\em SV-representation}.  In particular, it must be irreducible. 
\bigskip\\
\rem\label{rem: SV reps have maximal dimension}  Replacing $\rho$ by \textit{any} irreducible \rep  $\sigma:G\rightarrow \GL(W)$ in~(\ref{eq:1=s2/r+...}) gives $\dim(W)^2\leq |G/Z|$. This shows that SV-representations have the maximal possible dimension for irreducible representations of a given group $G$. If it has any such representations $G$  is said to be \textit{a group of central type}. For more details. 
and literature on such groups, see~\cite{Howe} and its references. 
\bigskip\\
Any representation $\rho:G\rightarrow \GL(V)$ is a composite 
$
G
\stackrel{\eta}{\rightarrow}G/\ker(\rho)\stackrel{\t \rho}{\rightarrow}\GL(V)
$
where $\t \rho$ is a faithful representation, and the surjectivity of the quotient map $\eta$ implies 
$
\eta(Z(G))\subset Z(G/\ker(\rho))
$. Thus we have a \textit{surjective} homomorphism, 
$$
\bar \eta : G/Z(G)=:\bar G \longrightarrow(G/(ker(\rho))/Z(G/\ker(\rho))=:\overline{G/(\ker(\rho))}
$$
and hence the equivalences 
\begin{eqnarray}
|\bar G|=|\overline{G/(\ker(\rho))}|&\Longleftrightarrow&
|\bar G|\leq |\overline{G/(\ker(\rho))}|\Longleftrightarrow
\mbox{$\bar \eta$ is an isomorphism} \Longleftrightarrow\nonumber\\
\mbox{$\bar \eta$ is injective} &\Longleftrightarrow&
\mbox{$Z(G)\supset\ker (\rho)$ and $Z(G)/\ker(\rho)=Z(G/\ker(\rho))$}
\label{eq:equiv condits re SV rho rho tilde}
\end{eqnarray}
Since $\rho$ is irreducible iff $\t \rho$ is irreducible, we have, by equivalent condition~\ref{part2:equivalent condits for SV} of 
Theorem~\ref{thm:equivalent condits for SV}:
\begin{lemma}\label{lemma:when rho and tilderho are SV}
	Let $\rho$, $\t \rho$, $\t \eta$, $\bar G$ and $\overline{G/(\ker(\rho))}$ be as above. If either of $\rho$ and $\t \rho$ is SV, then the other is SV iff any (hence all) of the equivalent conditions in~(\ref{eq:equiv condits re SV rho rho tilde}) holds.\ePf
\end{lemma}
Note that if $\rho$ is SV with central character $\psi$ then Theorem~\ref{thm:equivalent condits for SV}~\ref{part1:equivalent condits for SV} shows that
\begin{equation}\label{eq:ker of SV rep is ker of central char}
	\ker(\rho)=\ker(\psi)\subset Z
\end{equation}
since  $g\in G$ lies in $\ker(\rho)$ iff $\chi_\rho(g)=\dim(V)$. So one half of the last equivalent condition in~(\ref{eq:equiv condits re SV rho rho tilde}) is automatically satisfied in this case. 

Now suppose $\langle\cdot,\cdot\rangle$ is a  positive-definite, Hermitian form (or `PDHF') on $V$. We write $  {\rm U}(V)$ for the {\em unitary} subgroup of $\GL(V)$ w.r.t.\ this form (\ie\ those elements which preserve it) and call any representation $\rho:G\rightarrow \GL(V)$ `unitary' (w.r.t.\  $\langle\cdot,\cdot\rangle$) iff $\rho(g)\in   {\rm U}(V)$ for all $g\in G$. 
\begin{thm}[Isomorphism of SV-Representations]\label{thm:iso of gps gives iso of SV-R's}\ \\
Let $\rho_1:G_1\rightarrow \GL(V_1)$ and $\rho_1:G_1\rightarrow \GL(V_2)$ be two SV-representations 
with central characters $\psi_1$ and $\psi_2$ respectively, and $\alpha: G_1 \rightarrow G_2$ an isomorphism (so $\alpha(Z(G_1))=Z(G_2)$). Then 
\begin{enumerate}
	\item\label{part:Thm2_1} For an isomorphism $ T:V_1\rightarrow V_2$ satisfying 
	\begin{equation}\label{eq:T is iso of reps}
		T\circ\rho_1(g)=\rho_2(\alpha(g))\circ T
		\ \ \ \mbox{for all $g\in G_1$}
	\end{equation}
	to exist, it is necessary and sufficient that  
	\begin{equation}\label{eq: psi's related by theta}
		\psi_1=\psi_2\circ\alpha\ \ \ \mbox{on $Z(G_1)$}
	\end{equation} 
	\item\label{part:Thm2_2} If $T$ is an isomorphism satisfying~(\ref{eq:T is iso of reps}) then it is unique up to scalars.
	\item\label{part:Thm2_3} Suppose $\rho_1$ and $\rho_2$ are unitary
	w.r.t.\ PDHFs $\langle\cdot,\cdot\rangle_1$ on $V_1$ and $\langle\cdot,\cdot\rangle_2$ on $V_2$ respectively. If $T$ is an isomorphism satisfying~(\ref{eq:T is iso of reps}) than it can be taken to satisfy 
	\begin{equation}\label{eq:T preserves HIPs}
		\langle Tv, Tw\rangle_2=\langle v, w\rangle_1\ \ \ \mbox{for all $v,w\in V_1$}
	\end{equation}
	 (and is then unique up to multiplication by an element $t$ of $\bbT$, by~\ref{part:Thm2_2})
\end{enumerate}
\end{thm}   
\bPf\ The necessity of~(\ref{eq: psi's related by theta}) for~(\ref{eq:T is iso of reps}) is clear by taking $g\in Z(G_1)$. For sufficiency,  by condition~\ref{part1:equivalent condits for SV} of Theorem~\ref{thm:equivalent condits for SV}, equation~(\ref{eq: psi's related by theta})  implies that $\rho_1$ and $\rho_2\circ\alpha$ have the same character as representations of $G_1$, so they are isomorphic, \ie\ there exists an isomorphism $T$ satisfying~(\ref{eq:T is iso of reps}). Schur's lemma implies the uniqueness of $T$ up to $t\in\bbC^\times$. Finally, suppose the conditions of~\ref{part:Thm2_3} hold and $T$ is an isomorphism  satisfying~(\ref{eq:T is iso of reps}). We need to show that it can be adjusted by a scalar to satisfy~(\ref{eq:T preserves HIPs}). Let $T^\dagger$ denote the adjoint isomorphism $V_2\rightarrow V_1$ w.r.t.\ the two PDHFs. 
Since $\rho_1(g)$ and $\rho_2(\alpha(g))$ are unitary, taking adjoints everywhere in~(\ref{eq:T is iso of reps}) (w.r.t.\ the appropriate pairs of PDHFs) gives
$$
\rho_1(g)\inv\circ T^\dagger=T^\dagger\circ \rho_2(\alpha(g))\inv\ \ \ \mbox{for all $g\in G_1$}
$$  
and post-composing with~(\ref{eq:T is iso of reps}) on both sides, we get
$$
\rho_1(g)\inv\circ T^\dagger T\circ \rho_1(g)=T^\dagger T\ \ \ \mbox{for all $g\in G_1$}
$$
Now Schur's Lemma for $\rho_1$ implies $T^\dagger T=r\,{\rm id}_{V_1}\in\GL(V_1)$ for some  $r\in\bbC^\times$. But since 
$$r^\ast\langle v ,v\rangle_1=\langle T^\dagger T v , v\rangle_1=\langle T v ,T v\rangle_2$$ for any non-zero $v\in V_1$, we see that $r^\ast$, hence $r$, must be real and positive. Thus we may replace $T$ by $r^{-1/2}T$ to get $T^\dagger T=\id_{V_1}$ so 
$$
\langle Tv, Tw\rangle_2=\langle T^\dagger T v, w\rangle_1=\langle v, w\rangle_1\ \ \ \mbox{for all $v,w\in V_1$}
$$
as required.\ePf
\noindent Now we consider a fixed SV-representation $\rho:G\rightarrow   {\rm GL}(V)$ which is unitary w.r.t.\ a fixed PDHF $\langle\cdot,\cdot\rangle$ on $V$ and has central character $\psi$. First, we apply the Theorem with $\rho_1=\rho$ and  $\rho_2=\rho': G\rightarrow   {\GL}(V')$ another SV-representation  which is unitary w.r.t.\ some PDHF on $V'$, and  
$\alpha=\id_G$.  
We deduce that $\rho$ and $\rho'$ are  \textit{unitarily} equivalent iff they have the same central character. (In fact, this follows assuming only that $\rho'$ is irreducible, see~\cite[Prop 3.2.1]{Howe}.)

Next, we write $\Aut^{0}(G)$ for the subgroup of automorphisms in $\Aut(G)$ which \textit{fix $Z(G)$ elementwise} and take, $\rho_1=\rho_2=\rho$ and  $\alpha=\phi\in\Aut^{0}(G)$ in Theorem~\ref{thm:iso of gps gives iso of SV-R's}  to get 
\begin{equation}\label{eq:def T_phi}
	T_\phi\circ\rho(g)= \rho(\phi(g))\circ T_\phi
	\ \ \ \mbox{for all $g\in G$}
\end{equation}
for some $T_\phi\in   {\rm U}(V)$ determined up to an element of~\bbT. Given also $T_{\phi'}$ obeying a similar equation for another automorphism $\phi'\in \Aut^0(G)$ and applying it to both sides of~(\ref{eq:def T_phi}), we find, by uniqueness, that $T_{\phi'}T_{\phi}=T_{\phi'\phi}\bmod \bbT$. In other words, the map 
\begin{equation}\label{eq:def of Weil Rep}
	\displaymapdef{\cT_\rho}{\Aut^0(G)}{{\rm PU}(V)}{\phi}{T_\phi\bmod\bbT}
\end{equation}
to the projective unitary group is a well-defined homomorphism called the {\em Weil representation} associated to $\rho$. Note however that the above construction of $\cT_\rho$ is essentially \textit{inexplicit}.

In a certain sense, $\cT_\rho$ \textit{extends} $\rho$. Indeed, associated to any $h$ in $G$ is the inner automorphism $\phi_h$ in  $\Inn(G)\subset \Aut^0(G)$ sending every $g\in G$ to $hgh\inv$. 
When $\phi=\phi_h$, 
equation~(\ref{eq:def T_phi}) is clearly satisfied with $T_\phi=\rho(h)$. The map $\iota: G\rightarrow {\Aut^0}(G)$ sending $h$ to $\phi_h$,  has image ${\Inn}(G)$ and kernel $Z$. Thus we get a commuting diagram~(\ref{eq: first comm diag}) 
with exact rows
\begin{figure}[ht]
	\begin{equation}\label{eq: first comm diag}
		{\xymatrix{
			1\ar[r]&Z\ar[d]^{\psi}\ar[r]&G
			\ar[d]^{\rho}\ar[r]^{\iota}&
			{\rm Inn}(G)
			\ar[d]^{\cT_\rho}\ar[r]&1\\
			1\ar[r]&\bbT\ar[r]&{U}(V)\ar[r]&
			{\rm PU}(V)\ar[r]&1
		}}	
		\end{equation}
\end{figure}\ \\
\rem\label{rem: re the Clifford}\ Assume that $\psi$, hence $\rho$, is {\em injective} (\ie\ `faithful'). Then so is $\cT_\rho$. (Indeed, if $T_{\phi'}$ lies in $\bbT\id_V$  then~(\ref{eq:def T_phi}) -- with $\phi'$ for $\phi$ -- implies   $\phi'=\id_G$.) Now, (\ref{eq:def T_phi}) also shows that any $T_\phi$ lies in the normaliser $N:=N_{  {\rm U}(V)}(\rho(G))$. Conversely, our assumption implies easily that any $n\in N$ equals $T_\phi$ (defined modulo $\bbT$) for some automorphism $\phi\in\Aut^0(G)$. In summary, it implies that  $\cT_\rho$ maps $\Aut^0(G)$ isomorphically onto the image of $N$ in ${\rm PU}(V)$ and its subgroup ${\rm Inn}(G)$ onto the image of $\rho(G)$.   
\section{Overlap- and Angle-Maps, and Bouquets}\label{sec:Overlap- and Angle-Maps, and Bouquets}
\subsection{Definitions and First Properties} 
To start, we assume only that  $\rho:G\rightarrow \GL(V)$ is a complex representation of dimension $s<\infty$ of a finite group $G$ and that $\rho$ is unitary w.r.t.\ a PDHF $\langle\ ,\ \rangle$ on $V$. For each $u_1,u_2\in V$ we define a linear map 
$$
\displaymapdef{\cF_{u_1,u_2}}{V}{V}{u}{\sum_{g\in G} \langle  g\cdot u_1, u\rangle g\cdot u_2}
$$
\begin{prop}\label{prop:cF} For any $\rho$ as above and  $u_1,u_2\in V$
\begin{enumerate}
	\item \label{lemma:cF part1} $\cF_{u_1,u_2}(h\cdot u)=h\cdot\cF_{u_1,u_2}(u)$ for all $h\in G$ and $u\in V$
	\item \label{lemma:cF part2} $\Tr\,\cF_{u_1,u_2}=|G|\langle u_1,u_2\rangle$
\end{enumerate}
\end{prop}
\bPf\ Since $\rho$ is unitary,  
$\langle  g\cdot u_1, h\cdot u\rangle=\langle  h\inv g\cdot u_1, \cdot u\rangle$, and part~\ref{lemma:cF part1} follows easily. For part~\ref{lemma:cF part2}, let $\{e_1,\ldots,e_s\}$ be a an orthonormal basis of $V$. We find $\langle e_i, \cF_{u_1,u_2}(e_i)\rangle=\sum_{g\in G}
 \langle  g\cdot u_1, e_i\rangle \langle e_i,g\cdot u_2\rangle$ so that 
 \begin{eqnarray*}
 		\Tr\,\cF_{u_1,u_2}
&=&
 	\sum_{i=1}^s\sum_{g\in G}
 	\langle  g\cdot u_1, e_i\rangle \langle e_i,g\cdot u_2\rangle\\&=&
 	\sum_{g\in G}\sum_{i=1}^s
 	\langle  e_i,  g\cdot u_1\rangle^\ast \langle e_i,g\cdot u_2\rangle\\ &=&
 	\sum_{g\in G}\langle g\cdot u_1, g\cdot u_2\rangle\\ 
 	&=& |G|\langle u_1,u_2\rangle\ \ \ \Box 
 \end{eqnarray*}
 \begin{prop} Suppose $\rho$ as above is also {\em irreducible} and let $Z=Z(G)$. Then, for any $g\in G$ and $u_1,u_2,u\in V$, the vector $\langle  g\cdot u_1, u\rangle g\cdot u_2$ depends only on the the image $\bar{g}$ of $g$ in $G/Z$. Moreover, 
 \begin{equation}\label{eq: pre-pre-Clinometric}
 	\sum_{\bar{g}\in G/Z} \langle  g\cdot u_1, u\rangle g\cdot u_2=\left (\frac{|G/Z|\langle  u_1,u_2\rangle}{s}\right) u\addwords{for all $u_1,u_2,u\in V$}
 \end{equation}
 	\end{prop}
 \bPf\ Since $\rho$ is irreducible, $Z$ acts by roots of unity, which are of absolute value $1$, giving the first statement. For the second, Proposition~\ref{prop:cF} part~\ref{lemma:cF part1} and Schur's Lemma imply that $\cF_{u_1,u_2}$ must be multiplication by some $k\in \bbC$ and part~\ref{lemma:cF part2} implies $k=s\inv|G|\langle  u_1,u_2\rangle$. Thus
 $$
 |Z|\sum_{\bar{g}\in G/Z} \langle  g\cdot u_1, u\rangle g\cdot u_2=
 \cF_{u_1,u_2}(u)=\left (\frac{|G|\langle  u_1,u_2\rangle}{s}\right) u
 \ \ \ \ \ \ \ \ \ \ \ \ \ \ \Box 
 $$\smallskip\\
 \rem\ \textbf{(Tight Frames)}\label{rem:tight frames}\
Taking $u_1=u_2=v$, any non-zero vector in $V$, we find 
$$
\sum_{\bar{g}\in G/Z} \langle  g\cdot v, u\rangle g\cdot v=\left (\frac{|G/Z|\,|| v||^2}{s}\right) u\addwords{for all $u\in V$}
$$
In the language of~\cite{Waldron's book}, this means that, for any set $R$ of coset-representatives of $Z$ in $G$, the set $\{g\cdot v:g\in R\}$ is a  \textit{finite, equal-norm tight frame} in $V$. Such objects find applications in signal processing, quantum information theory, and elsewhere, see e.g.~\cite{Waldron's book}. 

 Note also that the rational number $|G/Z|/s$ appearing above is actually a positive \textit{integer}. Indeed, Theorem~3.12 of~\cite{Isaacs} states that, whenever $\rho$ is irreducible, $s$ divides the index in $G$ of the 
subgroup there denoted $Z(\chi)$, namely, the elements of $G$ acting by scalars through $\rho$. But irreducibility also implies  $Z\subset Z(\chi)$ because of the central character.\bigskip\\
%
For the rest of this section we shall continue to suppose that our representation $\rho:G\rightarrow V$ is \textit{irreducible} with central character $\psi:Z\rightarrow \mu_f$ where $Z=Z(G)$ and $f$ is the exponent of $Z$. We shall often abbreviate $G/Z$ to $\barG$ and write $\bar{g}$ for the image in $\barG$ of $g\in G$. 
We write $\bbP V$ for the projective space $\bbC^\times\backslash(V\setminus\{0\})$ identified with the set of all complex lines $\ell=\bbC v$ for some $v\neq 0$ in $V$.
Since $Z$ acts by $\psi$ on $V$, there is a well-defined action of $\barG$ on $\bbP V$ given  by $\bar{g}\cdot \ell:=\bbC (g\cdot v)$. For any $\ell\in\bbP V$ we write $\Gamma_\ell$ for the \textit{stabiliser} of $\ell$ in $\barG$, namely $\Gamma_\ell=\{\gamma\in\barG\,:\, \gamma\cdot\ell=\ell$\}. 

For any $\ell\in\bbP V$ and any map $\cR:\barG\rightarrow G$ which is right inverse to the quotient homomorphism $G\rightarrow \barG$, we define the {\em overlap-map}
\begin{equation}\label{eq:defn of Overlap-Map}
	\displaymapdef{\cAng_{\cR,\ell}}{\barG}{\bbC}{\gamma}{\langle v,\cR(\gamma)\cdot v\rangle}
\end{equation}
for any \textit{unit vector} $v\in\ell$. Note that $\cAng_{\cR,\ell}(\gamma)$ is independent of the choice of such $v$ and changing $\cR$ only multiplies it by $\psi(z)$ for some $z\in Z$ (so by an element of $\mu_f$).
We also define the \textit{angle} between any two lines $\ell,\ell'\in \bbP V$ to be 
$$
\fa(\ell,\ell')=|\langle v,v'\rangle|
$$
for any two unit vectors $v\in\ell$ and $v'\in \ell'$, so that $\fa(\ell,\ell')=\fa(\ell',\ell)$ and we define 
the {\em angle-map } $\fa_\ell$ by 
$$
\displaymapdef{\fa_{\ell}}{\barG}{\bbR_{\geq 0}}{\gamma}{\fa(\ell,\gamma\cdot \ell)
=|\cAng_{\cR,\ell}(\gamma)|\addwords{for any $\cR$}}
$$
\begin{prop}{\bf (Properties of the Overlap- and  Angle-Maps for an Irreducible Unitary Representation)}
\label{prop: properties of ov and ang}	\ \\
Let $\rho:G\rightarrow \GL(V)$ be irreducible and unitary of dimension $s$ and suppose $\ell$, $\ell'$ are in $\bbP V$. 
\begin{enumerate}
\item\label{part:properties of ov and ang0} We have 
 $0\leq \fa(\ell,\ell')\leq 1$ and $\fa(\ell,\ell')=1$ iff $\ell=\ell'$. In particular, the  angle-map $\fa_\ell$ takes values in $[0,1]$ and, for any $\gamma\in\barG$,  $\fa_\ell(\gamma)=1$ iff $\gamma\in \Gamma_\ell$. 
 \item\label{part:properties of ov and ang1} 
For every  $\gamma\in\barG$ we have 
\begin{equation}\label{eq: ov of inverse}
	\cAng_{\cR,\ell}(\gamma\inv)=
	\zeta\,\cAng_{\cR,\ell}(\gamma)^\ast
\end{equation}
for some $\zeta\in \mu_f$, 
so $\cAng_{\cR,\ell}(\gamma)\cAng_{\cR,\ell}(\gamma\inv)=\zeta\fa_{\ell}(\gamma)^2$ and 
$ 
\fa_{\ell}(\gamma\inv)=
\fa_{\ell}(\gamma)
$. 
\item\label{part:properties of ov and ang2.5} For all $\gamma,\gamma'\in \bar G$ we have 
\begin{equation}\label{eq: angle map on conjugates/image of line}
	\fa_{\gamma'\ell}(\gamma)=
	\fa_{\ell}({\gamma'}\inv \gamma \gamma')
\end{equation}
\item\label{part:properties of ov and ang2}
We have
\begin{equation}\label{eq:pre-Clinometric}
\sum_{\gamma\in\bar{G}}
\cAng_{\cR,\ell}(\gamma)^\ast\cAng_{\cR,\ell'}(\gamma)=
\frac{|G/Z|}{s}\fa(\ell,\ell')^2
\end{equation}
and \\
\begin{equation}\label{eq: The Clinometric Relation}
	\sum_{\gamma\in\bar{G}}
	\fa_{\ell}(\gamma)^2=
	\frac{|G/Z|}{s}
\end{equation}
\end{enumerate}
\end{prop}
\bPf The first statement in part~\ref{part:properties of ov and ang0} is a consequence of the complex Cauchy-Schwartz inequality and the others follow. For each $\gamma\in\barG$ we have 
 $\cR(\gamma\inv)=z\cR(\gamma)\inv$ for some $z\in Z$. Thus, for any unit vector $v$ generating $\ell$ we have: 
$$
\langle v,\cR(\gamma\inv)v\rangle=
\psi(z)\langle v,\cR(\gamma)\inv v\rangle=
\psi(z)\langle \cR(\gamma)v, v\rangle=
\psi(z)\langle v, \cR(\gamma)v\rangle^\ast
$$
by unitarity and hermitianness. This proves~(\ref{eq: ov of inverse}) and the rest of part~\ref{part:properties of ov and ang1} follows. Part~\ref{part:properties of ov and ang2.5} is clear from the definition and unitarity. 
For part~\ref{part:properties of ov and ang2}, take $u_1=u=v$ (a unit vector in $\ell$) and $u_2=v'$ (a unit vector in $\ell'$) in equation~(\ref{eq: pre-pre-Clinometric}) to get 
$$
\sum_{\gamma\in \barG} \langle  \cR(\gamma)\cdot v, v\rangle \cR(\gamma)\cdot v'=\left (\frac{|G/Z|\langle  v,v'\rangle}{s}\right) v
$$
and applying $\langle v',\ \rangle$ to both sides gives~(\ref{eq:pre-Clinometric}). Equation~(\ref{eq: The Clinometric Relation}) follows on taking $\ell'=\ell$.\ePf
\rem\label{rem: about inverses and inregrality}  Under  certain conditions (e.g.\ if $|\barG|$ is odd) we can choose $\cR$ such that   $\cR(\gamma\inv)=\cR(\gamma)\inv\ \forall\gamma\in\barG$ and the proof of~(\ref{eq: ov of inverse}) shows that we then have 
$
\cAng_{\cR,\ell}(\gamma\inv)=
\cAng_{\cR,\ell}(\gamma)^\ast\ \forall\gamma\in\barG
$.
\begin{defn}{\bf ($G$-Bouquets)}\label{defn:Bouquets}
	Let  $\rho:G\rightarrow \GL(V)$ be an irreducible unitary representation of dimension $s$ w.r.t.\ a PDHF $\langle\ ,\ \rangle$ on $V$.
	A {\em $G$-bouquet (with respect to $\rho)$} -- or just a {\em  bouquet} when $G$  and $\rho$ are clear -- is an orbit of $\bar G$ acting on $\bbP V$. We write $\PVmodbarG$ for the set of all such bouquets and say that a bouquet $\cY\in\PVmodbarG$ is {\em free} iff $\Gamma_\ell=\{1\}$ for one, hence any, $\ell\in\cY$, \ie\ iff $|\cY|=|\bar G|$. 
\end{defn}
\rem\ \textbf{(Polytope Picture of a Bouquet)}\label{rem:polytopes}\ 
We may imaginatively (and cautiously) associate the bouquet $\cY$  with a `complex polytope', say $\Poly(\cY)$, consisting of `vertices' (actually copies of $\bbT$) which are the intersections of the complex lines $\ell\in \cY$ with the $(2s-1)$-real-dimensional unit sphere in $V$ w.r.t.\ $\Iproddots$.  In this somewhat illusory  picture, $\bar{G}$ acts on the sphere by unitary  `rotations' which transitively permute the `vertices' of $\Poly(\cY)$. 
For $\ell\in\cY$ and $\gamma\in\bar G$, the `angle subtended at the origin' by the  `vertices' corresponding to the lines $\ell$ and $\gamma\cdot\ell$ is simply $\fa_\ell(\gamma)$.
\subsection{Equiangularity} 
\begin{defn}[Equiangularity of a Bouquet]\label{def:equiangularity} With $\rho$ as in Definition~\ref{defn:Bouquets}, a \emph{free} $G$-bouquet $\cY\in\PVmodbarG$ will be called {\em equiangular} iff 
$\fa(\ell,\ell')$ is independent of $\ell\neq \ell'\in\cY$. Equivalently, the angle map $\fa_\ell$ is constant on $\bar{G}\setminus\{1_{\barG}\}$ for one (hence any) $\ell\in\cY$. 
\end{defn}
\rem\ $G$ The condition  vacuously satisfied if $\bar G$ is trivial, \ie\ if $G$ is abelian, which is also equivalent to $s=1$ by irreducibility and freeness. In the Base Case detailed in~\S\ref{subsec:the Base Case} an equiangular bouquet is precisely a \textit{`Weyl-Heisenberg SIC(-POVM)'}\bigskip\\
Equiangularity is a very strong condition as the following result shows.
\begin{thm}\label{thm:equiangular existence and clinometric}
	Let $\rho$ be as in Definition~\ref{defn:Bouquets} and let $\cY$ be an equiangular $G$-bouquet for $\rho$. Then 
$\rho$ must be an SV-representation (see Section~\ref{reps of SV-type}). Moreover, it is necessary that  
	$\fa(\ell,\ell')=\frac{1}{\sqrt{s+1}}$ for all 
	$\ell\neq \ell'\in\cY$. 
\end{thm}
\bPf\ According to a theorem of~\cite{Delsarte et al}, the maximal number of equiangular lines in $V$ is $\dim(V)^2=s^2$. (Even when they are not, \textit{a priori}, the orbit of some group-action on $V$.) Since $\cY$ is free and equiangular, we deduce $|\barG|\leq s^2$. It follows that  $|\barG|=s^2$ (see Remark~\ref{rem: SV reps have maximal dimension}) and $\rho$ must be SV by Theorem~\ref{thm:equivalent condits for SV}~\ref{part2:equivalent condits for SV}. The second statement, holds vacuously if $s=1$. So assume $s>1$, let $a$ be the common value of $\fa(\ell_1,\ell_2)$ for 
$\ell_1\neq\ell_2\in\cY$ and choose $\ell\in\cY$. By freeness, equation~(\ref{eq: The Clinometric Relation}) reads  $1+(s^2-1)a^2=s$ so $a=(s+1)^{-\frac12}$. (In fact, this equality is also proven in~\cite{Delsarte et al} even without the group-orbit assumption.)\ePf
\section{Groups of Nilpotency Class $\leq 2$}\label{sec: nilpotent groups}
Throughout this section we suppose that 
\begin{equation}\label{eq:nilpotent class 2 condit}
\mbox{the group $G$ is nilpotent of class\,$\leq 2$}
\end{equation} 
This means that $\barG=G/Z$ is abelian. ($G$ is nilpotent of class 1 iff it is abelian.) As mentioned in the introduction, the Weyl-Heisenberg group ${\rm WH}(s)$ satisfies this condition, as indeed do all Generalised Heisenberg Groups, to be introduced in the next section. We shall see that it has significant simplifying consequences for overlap- and angle-maps and allows us to define `regular' bouquets. The SV condition -- and hence the existence of the Weil representation -- also becomes very natural
\subsection{SV Representations, Overlap- and Angle-Maps}   
Compare the following with~\cite[Theorem 13.5]{Waldron's book}.
\begin{thm}\label{thm:SV for class 2 nilpotent groups} 
If $G$ is nilpotent of class $\leq 2$ then every faithful, irreducible  representation $\rho: G \rightarrow \GL(V)$ is SV.
\end{thm}
\bPf\ This follows from Theorem~2.31 of~\cite{Isaacs} and  condition~\ref{part2:equivalent condits for SV} of our Theorem~\ref{thm:equivalent condits for SV}. (In fact, Isaacs' proof really establishes condition~\ref{partextra:equivalent condits for SV}.) It suffices to check that, in our situation, $Z$ is equal to the subgroup denoted $Z(\chi)$ in~\cite{Isaacs}, already mentioned in Remark~\ref{rem:tight frames}. But we noted there that irreducibility implies $Z\subset Z(\chi)$ and faithfulness clearly implies $Z\supset Z(\chi)$.\ePf
\rem\ If $\rho$ is faithful and irreducible, its central character must also be faithful, which forces $Z$ to be \textit{cyclic} (of order $f$).\bigskip\\
Applying Lemma~\ref{lemma:when rho and tilderho are SV} we deduce 
\begin{cor}
	If $G$ is nilpotent of class $\leq 2$ then an ireducible representation $\rho$ of $G$ is SV iff the equivalent conditions~(\ref{eq:equiv condits re SV rho rho tilde}) hold (with notations as in Lemma~~\ref{lemma:when rho and tilderho are SV}).\ePf 
\end{cor}
From now on, we suppose once again that the irreducible (but, for the moment, not necessarily SV) representation $\rho:G\rightarrow \GL(V)$   of dimension $s$ is unitary w.r.t.\ a PDHF $\langle\ ,\ \rangle$ on $V$, and write $\psi$ for its central character. Since $\bar G$ is abelian, equation~(\ref{eq: angle map on conjugates/image of line}) shows that, for any $\ell$ in $\bbP V$, the angle map $\fa_\ell$ depends only on the $G$-bouquet $\cY$ w.r.t.\ $\rho$ containing $\ell$ and may be denoted $\fa_\cY$. For the same reason we can write $\Gamma_\cY$ for the common stabiliser $\Gamma_{\ell}<\bar{G}$ of all $\ell\in\cY$. Now, condition~(\ref{eq:nilpotent class 2 condit}) is also  equivalent to the commutator $[g_1,g_2]=g_1g_2g_1\inv g_2\inv$ lying in $Z$ for all $g_1,g_2\in G$, so it descends to give a well-defined map of abelian groups: 
\begin{equation}\label{eq:def of deltaG for any G class 2 nilp}
	\displaymapdef{\delta_G}{\barG \times \barG}{Z}{(\barg_1,\barg_2)}{[g_1,g_2]}
\end{equation}
The standard commutator identities show that $\delta_G$ 
is skew-symmetric and bilinear. It is also non-degenerate in the sense that
 $\delta_G(\barg,\barg')=1_G$ for all $g'\in G$ clearly implies $g\in Z$, 
so $\barg=1_{\barG}$. For any unit vector $v\in\ell$, any $\cR$  and any $g,g'\in G$, we now have:   
\begin{eqnarray}
	\cAng_{\cR,\barg'\cdot \ell}(\barg)=
	\langle g'\cdot v,\cR(\barg)g'\cdot v\rangle&=&
	\langle g'\cdot v,[\cR(\barg),g']g'\cR(\barg)\cdot v\rangle\nonumber\\
	&=&
	\psi(\delta_G(\barg,\barg'))\langle g'\cdot v,g'\cR(\barg)\cdot v\rangle\nonumber\\
	&=&
	\psi(\delta_G(\barg,\barg'))\cAng_{\cR, \ell}(\barg)\label{eq:changing ell in overlaps under nilp}
\end{eqnarray}
A first consequence of~(\ref{eq:changing ell in overlaps under nilp}) is the existence of a map $\cAng_\cY$ (depending on $\ell$ only through $\cY$ and not at all on $\cR$):
\begin{equation}\label{eq:overlap map of bouquet}
	\displaymapdef{\cAng_\cY}{\barG}{\mu_f\backslash\bbC}{\barg}{\cAng_{\cR,\ell}(\barg) \bmod \mu_f=\langle v,g\cdot v\rangle \bmod \mu_f}
\end{equation}
We now introduce some notation for use in the rest of this paper. For any finite abelian group $X$ we shall write 
$\cM(X)$ for the  complex vector space of all functions from $X$ to $\bbC$
$$
\cM(X)=\cM(X,\bbC)=\{f:X\longrightarrow \bbC\}
$$
under pointwise addition and scalar multiplication, 
equipped with the PDHF given by 
\begin{equation}\label{eq:PDHF on M(A)}
	\langle f, g \rangle:=\sum_{x\in X}f(x)^\ast g(x)\addwords{for all  $f,g\in \cM(X)$}
\end{equation}  
Thus $\cM(X)$ has an obvious orthonormal basis $\cE_X:=\{\be_y\}_{y\in X}$ where 
$\be_y(x)$ is defined to be $1$ if $x=y$ and zero otherwise.  

Taking $X$ to be $\bar G$ as above, the central character $\psi$ of $\rho$ determines a map $\Upsilon_\psi\in\End_\bbC(\cM(\bar G))$: 
$$
\displaymapdef{\Upsilon_\psi}{\cM(\bar G)}{\cM(\bar G)}{f}
{\left(\gamma'\mapsto \sum_{\gamma\in \bar G }\psi(\delta_G(\gamma,\gamma'))f(\gamma)\right)}
$$
Thus, for any $\gamma_1,\gamma_2\in \bar G$, we have $\left\langle \be_{\gamma_2}, \Upsilon_\psi(\be_{\gamma_1})\right\rangle =\psi(\delta_G(\gamma_1,\gamma_2))$, so $\Upsilon_\psi$ is Hermitian by the skew symmetry of $\delta_G$. 
For the bouquet $\cY$ as above, the square $\fa_\cY^2=|\cAng_\cY|^2$ of the angle-map lies in $\cM(\bar G)$, is constant on cosets of $\Gamma_\cY$ in $\barG$ and, by Proposition~\ref{prop: properties of ov and ang}~\ref{part:properties of ov and ang0}, takes values in $[0,1]$ with
$\fa_\cY^2(\bar{\gamma})=1$ iff $\gamma\in\Gamma_\cY$.
\begin{thm}\label{thm:clinomteric thm for nilp gps}{\bf (Clinometry for Nilpotent Groups of Class $\leq 2$)}
With the above hypotheses and notations, $\fa_\cY^2$ is an eigenvector of $\Upsilon_\psi$ with eigenvalue $|\bar G|/s$, that is:
\begin{equation}\label{eq:aysquared is evec}
\Upsilon_\psi(\fa_\cY^2)=
(|\bar G|/s)\fa_\cY^2
\end{equation} 
\end{thm}
\bPf Choose any $\ell\in\cY$ and any $\cR$ right inverse to the quotient $G\rightarrow \barG$. Let $\gamma'$ be any element of $\barG$ and 
set $\ell'=\gamma'\cdot \ell$. Then equation~(\ref{eq:pre-Clinometric}) reads 
$$
\sum_{\gamma\in\bar{G}}
\cAng_{\cR,\ell}(\gamma)^\ast\cAng_{\cR,\gamma'\cdot \ell}(\gamma)=
(|\barG|/s)\fa_\cY(\gamma')^2
$$
But, according to~(\ref{eq:changing ell in overlaps under nilp}),  we have 
$
\cAng_{\cR,\gamma'\cdot \ell}(\gamma)=
\psi(\delta_G(\gamma,\gamma'))\cAng_{\cR,\ell}(\gamma) 
$, 
so 
$$
\sum_{\gamma\in\bar{G}}
\psi(\delta_G(\gamma,\gamma'))\fa_\cY(\gamma)^2=
(|\barG|/s)\fa_\cY(\gamma')^2
$$
for any $\gamma'\in\barG$, as required.\ePf 
\rem\ We shall refer to~(\ref{eq:aysquared is evec}) as the {\em  `clinometric relation'}. 
If we assume for simplicity that $\rho$ is faithful, then it confines the possibilities for the function $\fa_\cY^2$ to a subspace 
of $\cM(\barG)$ of a little over half its complex dimension. Indeed, in this case,  $\rho$ is SV by Theorem\ref{thm:SV for class 2 nilpotent groups} so $|\barG|=s^2$ and $\fa_\cY^2$ has eigenvalue $s$. On the other hand, since $\psi$ is injective and $\delta_G$ is non-degenerate, it is not hard to see that $\Upsilon_\psi=s^2 \id_{\cM(\barG)}$. (Check this on the $e_\gamma$'s for all $\gamma\in\barG$.) Thus $\cM(\barG)=\cM(\barG)^+\oplus \cM(\barG)^-$ where $\cM(\barG)^+$ and $\cM(\barG)^-$ are the eigenspaces of $\Upsilon_\psi$ with eigenvalues $s$ and $-s$ respectively. By considering also the trace of $\Upsilon_\psi$ we deduce easily that these eigenspaces have dimension $\half s(s+1)$ and $\half s(s-1)$ respectively. Note also that quation~(\ref{eq: The Clinometric Relation}) amounts (in the nilpotent case) to evaluating both sides of~(\ref{eq:aysquared is evec}) at $1_{\bar G}$. 
\subsection{Regularity and Invariance of a Bouquet}\label{subsec:Reg and Inv}
Any $\phi\in\Aut(G)$ preserves $Z$ and so descends to an automorphism  $\bar{\phi}$ of $\barG$. Obviously, $\Inn(G)$ actually fixes $Z$ elementwise, but  condition~(\ref{eq:nilpotent class 2 condit}) implies that its action on $\barG$ is also trivial. Thus $\bar{\phi}$ depends only on the image $\hat{\phi}$ of $\phi$ in $\Out(G)$. We use the resulting action of $\Out(G)$ on $\bar{G}$ to weaken the definition of equiangularity as follows.
\begin{defn}[$\cA$-Regularity of a Bouquet]\label{def:S-Regularity and Regularity}
	Let $\cY\in\PVmodbarG$ be a \emph{free} $G$-bouquet and $\cA$ be a subgroup of $\Out(G)$. We say that $\cY$ is {\em $\cA$-regular} (as in {\em `regular polytope'}) iff
	$\fa_\cY$ is constant on orbits of $\cA$ acting on $\barG$, that is, iff  
	\begin{equation}\label{eq:condition for S-regularity}
		\fa_\cY(\gamma)=\fa_\cY(\bar{\phi}(\gamma))\addwords{for all $\gamma\in \barG$ and $\hat{\phi}\in \cA$}
	\end{equation}
	In the case $\cA=\Out(G)$, we shall simply call $\cY$ {\em `regular'}.
\end{defn}
Since $\Out(G)$ acts on $\bar{G}$ by automorphisms, $\bar{G}\setminus\{1_{\bar{G}}\}$ is always a union of $\cA$-orbits, so equiangularity does imply $\cA$-regularity for any $\cA$. On the other hand, if $\bar{G}$
is not an elementary abelian $p$-group, $\bar{G}\setminus\{1_{\bar{G}}\}$ will contain elements of different orders, and so more than one orbit under $\Out^0(G)$. In such cases, it is \textit{a priori} possible to have regular bouquets which are not equiangular. 

For the rest of this section we shall add the condition that $\rho$ be SV (but not necessarily faithful). By Theorem~\ref{thm:iso of gps gives iso of SV-R's} and the comments following it, this implies that each $\phi\in\Aut^{0}(G)$ is associated to a unitary isomorphism $T_\phi:V\rightarrow V$ satisfying~(\ref{eq:def T_phi}), whose image in ${\rm PU}(V)$ is uniquely defined. In particular  $T_\phi$ acts unambiguously on $\bbP V$   and~(\ref{eq:def T_phi}) 
shows that it sends $G$-bouquets to $G$-bouquets. (In the polytope picture of Remark~\ref{rem:polytopes}, it corresponds to a `special rotation' about the origin that maps the polytope ${\rm Poly}(\cY)$ onto ${\rm Poly}(T_\phi(\cY))$ .)
\begin{prop}\label{prop:Upper Triangle Commutes}
Let $\rho:G\rightarrow\GL(V)$ be as above, let $\phi$ be an automorphism in $\Aut^0(G)$, choose $T_\phi$ as above and let $\cY\in\PVmodbarG$ be a $G$-bouquet.
\begin{enumerate}
\item\label{part:Upper Triangle Commutes 1} The bouquet $T_\phi(\cY)\in\PVmodbarG$ \emph{depends only on the image  $\hat{\phi}$ of $\phi$ in $\Out^0(G)$} and 
$$
\displaynomapdef{\Out^{0}(G)\times\PVmodbarG}{\PVmodbarG}{(\hat{\phi}\, ,\, \cY)}{T_\phi(\cY)}
$$ 
is an action of $\Out^0(G)$ on the set $\PVmodbarG$.
\item Let $\cR:\barG\rightarrow G$ be any right inverse to the quotient $G\rightarrow \barG$. 
Then $\cR^\phi:=\phi\circ\cR\circ\bar{\phi}\inv$ 
is also such a right inverse, and for all $\ell\in\bbP V$ the following diagram commutes
\begin{equation}\label{diag:Upper Triangle Commutes}
{\xymatrix{
	\bar{G}\ar[dd]^{\bar{\phi}}\ar[rrrr]^{\cAng_{\cR,\ell}}&&&&\bbC\\
	&&&&\\
	\bar{G}
	\ar[rrrruu]_{\ \ \cAng_{\cR^\phi, T_\phi(\ell)}}&&&&
}}
\end{equation}
\item\label{part:Upper Triangle Commutes 3} The maps $\cAng_{T_\phi(\cY)}\circ\bar{\phi}$ and $\cAng_{\cY}$ from $\bar{G}$ to $\mu_f\backslash\bbC$ are the same, as are the maps  $\fa_{T_\phi(\cY)}\circ\bar{\phi}$ and $\fa_{\cY}$ from $\bar{G}$ to $[0,1]$. In particular, the maps $\cAng_{T_\phi(\cY)}$ and $\cAng_{\cY}$ have the same image, as do the maps $\fa_{T_\phi(\cY)}$ and $\fa_{\cY}$. 
\end{enumerate}
\end{prop}
\bPf For every inner automorphism $\phi_h\in\Inn(G)$ , we have seen (\cf\ diagram~(\ref{eq: first comm diag})) that $T_{\phi_h}=\rho(h)$ up to scalars, so $T_{\phi_h}(\cY)=\cY$, since $\cY$ is an orbit of $\bar{G}$. The first part follows since $\cT_\rho$ is a homomorphism.
For the second, it is clear that $\cR^\phi$ has the required property. Let $v$ be a unit vector generating $\ell$ so that $T_\phi(v)$ generates $T_\phi(\ell)$. Since $T_\phi$ is unitary, $T_\phi(v)$ is also a unit vector and, using equation~(\ref{eq:def T_phi}) we have, for any $\gamma\in \barG$, 
$$
\cAng_{\cR^\phi, T_\phi(\ell)}(\bar\phi(\gamma))=
\langle T_\phi(v),\phi(\cR(\gamma))\cdot T_\phi(v) \rangle=
\langle T_\phi(v), T_\phi(\cR(\gamma)\cdot v) \rangle=
\langle v, \cR(\gamma)\cdot v\rangle=\cAng_{\cR, \ell}(\gamma)
$$ 
giving~(\ref{diag:Upper Triangle Commutes}). Taking $\ell$ to be any line in the bouquet $\cY$, the third part also follows from the above equation and the definitions of the overlap- and angle-maps of $\cY$.\ePf 
\noindent We can now reformulate the regularity of a bouquet $\cY$ with respect to a subgroup $\cA$ of $\Out^0(G)$ in terms of the equality of the angle maps associated to the image-bouquets $T_\phi(\cY)$ for the corresponding automorphisms
$\phi\in\Aut^0(G)$. Indeed, parts~\ref{part:Upper Triangle Commutes 1} and~\ref{part:Upper Triangle Commutes 3} of Proposition~\ref{prop:Upper Triangle Commutes} lead immediately to the following. 
\begin{prop}\label{prop:equiv condit for regularity}{\bf (Equivalent Condition for $\cA$-Regularity in the SV Case)}\ \\
	Let $\cY\in\PVmodbarG$ be a \emph{free} $G$-bouquet and $\cA$ be a subgroup of $\Out^0(G)$. Then $\cY$ is {\em $\cA$-regular} iff 
	$$
	\fa_{T_\phi(\cY)}=\fa_{\cY}:\barG\rightarrow[0,1]
	\addwords{for all $\hat\phi\in\cA$}
	$$
	where, as usual, $\phi$ denotes any lift of $\hat\phi$ to $\Aut^0(G)$ and the bouquet $T_\phi(\cY)$ is independent of the lift chosen.\ePf  
\end{prop}
In practice, the isomorphisms $T_\phi$ may be harder to calculate than the action of $\cA$ on $\bar{G}$. In such cases, formulation~(\ref{eq:condition for S-regularity}) of $\cA$-regularity, which is intrinsic to $\cY$, may be more useful than that of Proposition~\ref{prop:equiv condit for regularity}. Next, we make the 
\begin{defn}{\bf (Symmetry Group and Invariance of a Bouquet)}\\
Let  $\cY\in\PVmodbarG$ be a $G$-bouquet. The \emph{symmetry group} of $\cY$ is its stabiliser in $\Out^0(G)$ for the action of part~\ref{part:Upper Triangle Commutes 1} of Proposition~\ref{prop:Upper Triangle Commutes} namely the subgroup 
$$
{\rm Sym}(\cY):=\{\hat{\phi}\in\Out^0(G)\,:\, T_\phi(\cY)=\cY\} < \Out^0(G)
$$
Let $\hat \cE$ be a subrgroup of $\Out^0(G)$. We say that $\cY$ is 
\emph{$\hat \cE$-invariant} iff $T_\eps(\cY)=\cY$ for all $\eps\in\Aut^0(G)$ with $\hat\eps\in \hat \cE$, that is, iff  $\hat \cE\subset {\rm Sym}(\cY)$.
\end{defn}
Proposition~\ref{prop:equiv condit for regularity} shows that the $\hat \cE$-invariance of a free bouquet implies its $\hat \cE$-regularity. In fact, something stronger follows immediately from Proposition~\ref{prop:Upper Triangle Commutes}\ref{part:Upper Triangle Commutes 3}: 
\begin{prop}
If $\cY$ is $\hat \cE$-invariant then $\frak O_\cY:\bar G\rightarrow \mu_f\backslash\bbC$ is constant on the orbits of $\hat \cE$ acting on $\bar G$.$\ \Box$ 
\end{prop}
\rem\ 
In the Base Case, the symmetry groups of {\em all  known equiangular ${\rm WH}(s)$-bouquets
 with $s>3$}  (see \S\ref{subsec:the Base Case}) are non-trivial. Indeed, it has been observed experimentally that they all contain an an element $\widehat{\phi_\fz}$ of order $3$. Zauner conjectured that this happens for all such bouquets, and specified a matrix of order $3$ that should represent a corresponding $T_{\phi_\fz}$ in a standard basis of $V$. See his PhD thesis~\cite{Zauner}.
\section{Generalised Heisenberg Groups and Their Schr\"odinger Representations}\label{sec:GHGs and their SRs}
In this section we introduce a large class of abstract, class $\leq 2$ nilpotent groups with explicit, irreducible, unitary, SV-representations, for which we can apply the preceding theory. In Section~\ref{sec:Examples} we shall provide a broad family of concrete `arithmetic' examples of such groups.
\subsection{Definitions and Basic Properties}
Given any three \textit{finite} abelian groups 
 $A$, $B$ and $C$, and a ($\bbZ$-)bilinear map $\lambda:A\times B\rightarrow C$ we define the {\em generalised Heisenberg group} -- or \textit{GHG} for short -- $\cH=\cH(A,B,C,\lambda)$ attached to the tuple $(A,B,C,\lambda)$ to be the set $A\times B\times C$ with group-law   
$$
h(a,b,c)h(a',b',c')=h(a+a', b+b',c+c'+\lambda(a,b') )
$$ 
(The symbol `$h$' is redundant but helpful.) The following are easily checked: $\cH$ is a finite group with $h(a,b,c)\inv=h(-a,-b,\lambda(a,b)-c)$ and
\begin{equation}\label{eq:generators of HABC}
	h(a,b,c)=h(0,0,c)h(0,b,0)h(a,0,0)=h(0,b,0)h(a,0,0)h(0,0,c)
\end{equation}
The natural surjection $\pi:\cH\rightarrow A\oplus{}B$ and injection $m:C\rightarrow \cH$ (with image $h(0,0,C)$)  are homomorphisms. If $h=h(a,b,c)$ and 
$h'=h(a',b',c')$ then 
\begin{equation}\label{eq: commutators in HeiStan}
	[h,h']=h(0,0,\lambda(a,b')-\lambda(a',b))
\end{equation}
The centre $Z$ and derived subgroup $\cH'$ of $\cH$ are therefore given by 
\begin{equation}\label{eq:centre and commutator}
	Z=h(K_A, K_B, C)\supset m(C) \supset  
	\cH'= h(0,0,\Lambda) 
\end{equation}
where $K_A$ and $K_B$ denote the left- and right-kernels of $\lambda$, and $\Lambda$ denotes the image of the map 
$A\otimes{}B\rightarrow C$ defined by $\lambda$. In particular, $\cH$ is always \textit{nilpotent of class $2$} (or $1$ iff $\lambda=0$) and we have a \textit{central} extension of groups
\begin{equation}\label{eq:ses of gps}
	1\rightarrow C\stackrel{m}{\longrightarrow} \cH\stackrel{\pi}
	\longrightarrow 
	A\oplus{}B\rightarrow 1
\end{equation}
There are canonical, commuting automorphisms $\phi_-$, $\phi^-$ and $\phi_{-1}=\phi^-\circ\phi_-$ in $\Aut(\cH)$:
$$
\phi_-(h(a,b,c))=h(-a,b,-c),\ \  
\phi^-(h(a,b,c))=h(a,-b,-c)\ \ \mbox{and}\ \   
\phi_{-1}(h(a,b,c))=h(-a,-b,c)
$$ 
all of order $1$ or $2$. 
Thus ${\frak V}:=\{\id_\cH, \phi_-, \phi^-, \phi_{-1}\}$ is a subgroup of $\Aut(\cH)$ isomorphic to (a quotient of)  $(\Zover{2})^2$. 
Defining $\lambda\op: B\times A \rightarrow C$
by $\lambda\op(b,a)=\lambda(a,b)$ we get an \textit{anti}-isomorphism $\t \varphi:
\HeiStan\rightarrow\Hei{B}{A}{C}{\lambda\op}$ sending $h(a,b,c)$ to $h(b,a,c)$. Composing $\t \varphi$ with inversion therefore gives an isomorphism 
$$
\displaymapdef
{\varphi}
{\HeiStan}
{\Hei{B}{A}{C}{\lambda\op}}
{h(a,b,c)}
{h(-b,-a,\lambda(a,b)-c)}
$$
In special cases where $A=B$ and $\lambda$ is symmetric (\eg\ the Base Case for odd $d$) $\varphi$ becomes an additional automorphism of $\Hei{A}{A}{C}{\lambda}$ of order $1$ or $2$. 
Together with $\frak V$ it then generates a subgroup of $\Aut(\Hei{A}{A}{C}{\lambda})$  isomorphic to (a quotient of) the dihedral group of order $8$ and with $\phi_{-1}$ in its centre and in $\Aut^0(\Hei{A}{A}{C}{\lambda})$.\bigskip\\
\rem\ In another wide class of cases (Hypothesis~ND-C) there are also {\em non-canonical} isomorphisms between $A$ and $B$ and hence between ${\HeiStan}$ and ${\Hei{B}{A}{C'}{\mu}}$ for any {\em nondegenerate} $\mu$  and $C\cong C'$ (see Corollary~\ref{cor:implication about iso} with $R=\bbZ$).
\bigskip\\
Let $p:C\rightarrow \bbC^\times$ be any homomorphism. Since $C$ is finite, its values are roots of unity in $\bbT$. For any  
$f\in\cM(A)$ (see \S\ref{sec: nilpotent groups}) and any $h(a,b,c)\in\cH$ we define 
$h(a,b,c)\cdot_{\sigma_p}f\in\cM(A)$ by  
\begin{equation}\label{eq: defn of Schrod psi-action on functions}
	\left(h(a,b,c)\cdot_{\sigma_p}f\right)(x)= 
	p(\lambda(x,b)+c))f(x+a)\ \ \ \mbox{for all $x\in A$}
\end{equation}
One checks easily that `$\cdot_{\sigma_p}$' thus defined is a $\bbC$-linear, left group-action of  $\cH$ on $\cM(A)$ and that it preserves its PDHF. 
It follows that we get a well defined a unitary representation
$$
\sigma_p:\cH=\HeiStan\longrightarrow \U(\cM(A))
$$ 
by setting $\sigma_p(h)(f):=h\cdot_{\sigma_p} f$ for any $h\in\cH$ and $f\in\cM(A)$. We call this the {\em (discrete) left Schr\"odinger \rep} of $\cH$, attached to $p$. 
One checks similarly that there is another $\bbC$-linear left 
group action `$\cdot_{\tau_p}$' of $\cH$ on $\cM(B)$ given for any $l\in\cM(B)$ by 
\begin{equation}\label{eq: RIGHT Schrod psi-action on functions}
	\left(h(a,b,c)\cdot_{\tau_p}l\right)(y)= 
	p(c-\lambda(a,y+b)))l(y+b)\ \ \ \mbox{for all $y\in B$}
\end{equation}
and that it preserves the PDHF on $\cM(B)$ defined in the same way as~(\ref{eq:PDHF on M(A)}). 
Thus we get another well-defined a unitary representation
$$
\tau_p:\cH\longrightarrow \U(\cM(B))
$$ 
by setting $\tau_p(h)(l):= h\cdot_{\tau_p}l$ for any $h\in\cH$ and $l\in\cM(B)$. We call this the {\em (discrete) right Schr\"odinger \rep} of $\HeiStan$ attached to $p$. 
Note that the central element $m(c)=h(0,0,c)$ of $\cH$ acts by multiplication by $p(c)$ in both Schr\"odinger representations.\bigskip\\ 
\rem\label{rem:linking tau to sigma of HBA gp via varphi}\ In fact, it is an easy exercise to check that $\tau_p=\t \sigma_{p^\ast}\circ\varphi\circ\phi_{-1}$ where where $\t \sigma_{p^\ast}$ denotes the left Schr\"odinger representation of $\Hei{B}{A}{C}{\lambda\op}$ with respect to the complex conjugate homomorphism $p^\ast$, \ie\ $p$ composed with inversion in $C$.\bigskip\\
Calculating traces in the bases $\cE_A$ and $\cE_B$, one shows easily:  
\begin{lemma}\label{lemma:char of rho psi} 
	The characters
	$\chi_{\sigma_p}$ and $\chi_{\tau_p}$ and  of $\sigma_p$ and $\tau_p$ respectively are given by 
	$$
	\chi_{\sigma_p}(h(a,b,c))=
	\eitherortw{p(c)|A|}{if $a=0$ and $\lambda(A,b)\subset\ker(p)$}{0}{otherwise}
	$$
	and
	$$
	\chi_{\tau_p}(h(a,b,c))=
	\eitherortw{p(c)|B|}{if $b=0$ and $\lambda(a,B)\subset\ker(p)$}{0}{otherwise}
	$$
	for all $h(a,b,c)\in\cH$.\ePf
\end{lemma}
For reasons that will become apparent in the next subsection, we will mainly be concerned with $\sigma_p$ in the rest of this one.
In fact, we shall be concerned more generally with homomorphisms, isomorphisms and automorphisms of various kinds between generalised Heisenberg groups, and how they interact with $\sigma_p$. A particularly simple type is as follows. Let $(A,B,C,\lambda)$ and 
$({\t A},{\t B},{\t C},{\t \lambda})$ be two tuples as above, suppose $t_A:A\rightarrow {\t A}$, $t_B:B\rightarrow {\t B}$, $t_C:C\rightarrow {\t C}$ are homomorphisms and consider the map 
$$
\displaymapdef{t_A\times t_B\times t_C}{\HeiStan}{\cH({\t A},{\t B},{\t C},{\t \lambda})}
{h(a,b,c)}{h(t_A(a),t_B(b),t_C(c))}
$$
It is easy to prove
\begin{prop}\label{prop: criterion for diagonal homs}{\bf (Criterion for Diagonal Homomorphisms)}
	With notations and assumptions as above, 
	$t_A\times t_B\times t_C$ is a homomorphism if and only if 
	\begin{equation}\label{eq:compatibility of lambdas}
		t_C(\lambda(a,b))={\t \lambda}(t_A(a), t_B(b))\addwords{for all $a\in A$, $b\in B$\ \ \ \ \  \ \ $\Box$}
	\end{equation} 
\end{prop}
Homomorphisms of form $t_A\times t_B\times t_C$ will be called \textit{diagonal}. They commute  with left Schr\"odinger representations as follows.  Consider the linear map 
$\check t_A:\cM(A)\rightarrow \cM({\t A})$ sending $\be_a$ to $\be_{t_A(a)}\ \forall a\in A$. Explicitly 
\begin{equation}\label{eq:distribution relation}
	\check t_A(f)({\t a}):=\sum_{a\in A \atop  t_A(a)={\t a}}f(a)\addwords{for all ${\t a}\in {\t A}$}
\end{equation}
where we interpret the R.H.S.\ as zero for ${\t a}\nin t_A(A)$.\bigskip\\
\rem\label{rem:injectivity of tAstar} It is easy to see that $\check t_A$ preserves PDHFs iff $t_A$ is injective.
\begin{prop}\label{prop:"for future use"} In the situation of Proposition~\ref{prop: criterion for diagonal homs}, suppose $\underline{t}:=t_A\times t_B\times t_C$ is a diagonal homomorphism (so~(\ref{eq:compatibility of lambdas}) holds) and suppose $p:C\rightarrow\bbC^\times$ and ${\t p}:{\t C}\rightarrow\bbC^\times$ are homomorphisms satisfying ${\t p}\circ t_C=p$. Then  
	\begin{equation}\label{eq:actions commute with projection}
		\underline{t}(h)\cdot_{\sigma_{\t p}} \check t_A(f) = \check t_A (h\cdot_{\sigma_p} f)
	\end{equation} 
for all $h\in\HeiStan$ and $f\in\cM(A)$. 
\end{prop} 
\bPf\ By linearity, it suffices to check~(\ref{eq:actions commute with projection}) for $f$ every basis element $\be_a$  for $a\in A$. Setting $h=h(a_1,b_1,c_1)$ and evaluating both sides at $x'\in {\t A}$, we find they both vanish unless $x'=t_A(a-a_1)$ in which case they both equal 
$p(\lambda(a-a_1,b_1)+c_1)$. \ePf
\rem\label{rem:p-adic measures}\ {\bf (Projective Systems of GHGs)}\ Notice that~(\ref{eq:distribution relation}) is a form of distribution relation. One can also  study certain profinite Heisenberg groups as projective limits of systems of the form $(\cH_i:=\Hei{A_i}{B_i}{C_i}{\lambda_i})_{i\in\cI}$, with diagonal transition homomorphisms to get 
$$
\lim_{\leftarrow}\Hei{A_i}{B_i}{C_i}{\lambda_i}\cong \Hei{A_\infty}{B_\infty}{C_\infty}{\lambda_\infty}
$$
where $A_\infty:=\displaystyle{\lim_\leftarrow A_i}$ etc. This is somewhat in line with the adelic Heisenberg groups of~\cite{Weil}. In such situations, the distribution relations mean that the (infinite-dimensional)  spaces for the left Schr\"odinger representations of such profinite groups can be naturally interpreted as   \textit{distributions, or measures,} contained in 
$\cM_\infty:={\displaystyle \lim_{\leftarrow}}\,\cM(A_i)$ instead of the real- or complex-valued {\em functions}, on the compact set  $A_\infty$, which appear in the situations studied in~\cite{Weil}. A $p$-adic example of this set-up is outlined in subsection~\ref{subsec:perspectives}. 
\subsection{The `Non-Degenerate-by-Cyclic' Condition}\label{subsec:ND-C}
Unless otherwise stated we assume henceforth the `ND-C' condition on the tuple $(A,B,C,\lambda)$ which is commonly satisfied in the applications.  
\begin{hyp}[The Condition ND-C]\label{hyp1}
	\begin{enumerate}\item[] 
	    \item\label{hyp1:three} $\lambda$ is non-degenerate (\ie\ $K_A=K_B=0$)
		\item\label{hyp1:two} $C$ is cyclic, say $C\cong \Zover{r}$ with $r\in\bbZ_{>0}$. 
		\end{enumerate}
\end{hyp} 
Since $C$ is cyclic, we can choose an injective homomorphism $q:C\rightarrow\bbC^\times$ to get \begin{equation}\label{eq:dual maps}
A\stackrel{\lambda_\ast}{\longrightarrow}\Hom(B,C)\stackrel{q_\ast}{\longrightarrow}\Hom(B,\bbC^\times ):=\hat{B}
\end{equation}
where the dual group $\hat{B}$ is non-canonically isomorphic to $B$ and the homomorphism $\lambda_\ast$ (induced by $\lambda$) is injective by~\ref{hyp1:three}, as is $q_\ast$. We deduce $|A|\leq |\hat{B}|=|B|$. The same argument but reversing the roles of $A$ and $B$  gives $|B|\leq |A|$. Thus $|A|=|B|$ so $\lambda_\ast$ and $q_\ast$ must be  \textit{isomorphism}s and $A\cong\hat{B}\cong B$. 
We write $s$ and $e$ respectively for the common cardinalities and exponents of $A$ and $B$. Thus $s$ is a multiple of $e$ and has the same prime factors. Since $r$ kills $\Hom(B,C)\cong\hat{B}$, it must also be a multiple of $e$. If $\lambda$ is surjective then $e$ will kill $C$, so also $r|e$ and $r=e$. Note that we \textit{do not} insist that this be the case under condition ND-C. Thus, in particular, it is quite possible for $r$ to have prime factors that do not divide $e$.

Abbreviating $\HeiStan$ to $\cH$ as before,   
we have $Z=Z(\cH)=m(C)$ by~(\ref{eq:centre and commutator})  and $\cH'$ is the unique subgroup of $Z$ of cardinality $e$. Also,
\begin{equation}\label{eq:all set up for an SV rep}
\dim(\cM(A))=s\ \ \ \mbox{and}\ \ \ |\cH:Z|=s^2
\end{equation}
since $\cH/Z\cong A\oplus B$. Clearly, $\cH$ has cardinality $s^2 r$ and we shall see in \S\ref{subsec:Aut Gps via D} that, at least when $r$ is odd, it is precisely the exponent of $\cH$.

By Theorem~\ref{thm:SV for class 2 nilpotent groups}, faithful irreducible representations and faithful SV representations of $\cH$ are the same thing. We can now characterise them precisely up to isomorphism.
\begin{prop}\label{prop:faithful SV reps for HeiStan under Hyp 5.1}
Suppose condition ND-C is satisfied.
\begin{enumerate} 
\item\label{part:faithful SV reps for HeiStan under Hyp 5.1 1} 
If $p:C\rightarrow \bbC^\times$ is an injective homomorphism then $\sigma_p$ is a faithful SV representation of $\cH$ with central character $\psi$ given by $\psi\circ m=p$. 
\item\label{part:faithful SV reps for HeiStan under Hyp 5.1 2} 
If $\rho$ is any faithful irreducible representation of $\cH$ with central character $\psi$, then $\psi$ is injective and $\rho$ is isomorphic to $\sigma_p$ iff $p=\psi\circ m$. 
\end{enumerate}
\end{prop}
\bPf\ If $p$ is injective then, by the non-degeneracy of $\lambda$ and Lemma~\ref{lemma:char of rho psi},  $\chi_{\sigma_p}(h(a,b,c))=0$ unless $a=b=0$, \ie\ $h(a,b,c)\in m(C)=Z$. Thus $\sigma_p$ is SV by equation~(\ref{eq:all set up for an SV rep}) and condition~\ref{part1:equivalent condits for SV} of Theorem~\ref{thm:equivalent condits for SV}. It is faithful by~(\ref{eq:ker of SV rep is ker of central char}) and $\psi\circ m=p$ also follows from Lemma~\ref{lemma:char of rho psi}. This proves part~\ref{part:faithful SV reps for HeiStan under Hyp 5.1 1}. For part~\ref{part:faithful SV reps for HeiStan under Hyp 5.1 2}, $\psi$ is obviously injective and Theorem~\ref{thm:SV for class 2 nilpotent groups} tells us $\rho$ is SV. It will be isomorphic to $\sigma_p$ iff it has the same character, which is equivalent to $p=\psi\circ m$ by Theorem~\ref{thm:equivalent condits for SV}\ref{partextra:equivalent condits for SV}.\ePf 
\noindent If ND-C holds and $p$ is injective, then Lemma~\ref{lemma:char of rho psi}  shows that $\tau_p$ has the same character as $\sigma_p$, so the two are abstractly isomorphic. We now give an explicit unitary isomorphism between them under these conditions. Indeed, we define the $\bbC$-linear map $\xi_p:\cM(A)\rightarrow\cM(B)$:  
$$
\xi_p(f)(y)=s^{-\frac12}\sum_{x\in A} f(x)p(\lambda(x,y))
\addwords{for all $y\in B$, for all $f\in \cM(A)$}
$$
Thus $\xi_p$ is the normalised, discrete Fourier transform. Given $a\in A$,  define $\omega_{a,p}\in \hat B$ and $\t \omega_{a,p}\in \cM(B)$ by
$$
\omega_{a,p}(y):=p(\lambda(a,y))\addwords{for all $y\in B$,  and}\ \ \ \t \omega_{a,p}:=s^{-\frac12}\omega_{a,p}
$$ 
Note that $\omega_{a,p}$ is just the image of $a$ under the isomorphisms of~(\ref{eq:dual maps}) (with $p$ for $q$) so it runs through $\hat B$ as $a$ runs though $A$. Hence, by the character theory of finite abelian groups, it follows that $\cG_B:=\{\t \omega_{a,p}:a\in A\}$ is an orthonormal basis of $\cM(B)$ independent of $p$. Similarly, we define 
$\omega_{b,p}\in \hat A$ for each $b$ by $\omega_{b,p}(x):=p(\lambda(x,b))$, set $\t \omega_{b,p}=s^{-\frac12}\omega_{b,p}\in \cM(A)$ and get an orthonormal basis $\cG_A:=\{\t \omega_{b,p}:b\in B\}$ of $\cM(A)$. 
\begin{prop}\label{prop:isomorphism of left and right Schrod reps} {\bf($\tau_p$ and $\sigma_p$ are Fourier Duals under Condition ND-C)}\\
	With $p$ and $\xi_p : \cM(A)\longrightarrow\cM(B)$ as above, 
	\begin{enumerate}
		\item\label{part:isomorphism of left and right Schrod reps1} 
		$\xi_p(e_a)=\t \omega_{a,p}\ \forall a\in A$ and $\xi_p(\t \omega_{p,b})=e_{-b}\ \forall b\in B$ and 						 
		\item\label{part:isomorphism of left and right Schrod reps2}
		 $\xi_p$ is an isomorphism satisfying 
		$\langle \xi_p(f_1), \xi_p(f_2)\rangle=\langle f_1, f_2\rangle\  \forall f_1,f_2\in \cM(A)$ 
		\item\label{part:isomorphism of left and right Schrod reps3} We have
		\begin{equation}\label{eq:xip is an iso of schrod reps}
			\xi_p(h\cdot_{\sigma_p}f)= h\cdot_{\tau_p}\xi_p(f)
			\addwords{for all $h\in\HeiStan$ and $f\in \cM(A)$}
		\end{equation}
	\end{enumerate}
\end{prop}
\bPf\ The formulae in~\ref{part:isomorphism of left and right Schrod reps1} are easily checked and imply~\ref{part:isomorphism of left and right Schrod reps2} since, for example, $\xi_p$ takes $\cE_A$ to $\cG_B$ and both are orthonormal bases. For part~\ref{part:isomorphism of left and right Schrod reps3}, since `$\cdot_{\sigma_p}$', and `$\cdot_{\tau_p}$' are both group actions,   
it suffices to check~(\ref{eq:xip is an iso of schrod reps}) taking $h$ to be each of the generators of $\HeiStan$ of form $h(0,0,c)$, $h(0,b,0)$ and $h(a,0,0)$ 
(see~(\ref{eq:generators of HABC})). This is straightforward and LTR.
\ePf
\noindent It follows that $\bbP\xi_p$ maps every $\cH$-bouquet  $\cY\subset\bbP\cM(A)$ w.r.t.\ $\sigma_p$ explicitly onto a `dual' $\cH$-bouquet $\bbP\xi_p(\cY)\subset\bbP\cM(B)$ w.r.t.\ $\tau_p$ with the same overlap- and angle-maps. For the purposes of studying equiangular and $\cA$-regular bouquets, $\sigma_p$ and $\tau_p$ are are therefore equivalent so we focus on $\sigma_p$ from now on, often calling it simply \textit{the} Schr\"odinger representation of $\cH$. 
\subsection{GHGs `with $R$-Structure'}
We temporarily suspend the assumption of condition ND-C to greatly enrich the theory of generalised Heisenberg groups: 
we consider tuples $(A,B,C,\lambda)$  
where\textit{ $A$ and $B$ are modules over a commutative ring $R$}. (Note: we allow $R$ to be infinite, but its action on $A$ will factor through the necessarily finite ring $R/\ann(A)$.)
We also require $\lambda$ to be {\em $R$-balanced} in the sense that 
\begin{equation}\label{eq: lambda is R-balanced}
	\lambda(ra,b)=\lambda(a,rb)\addwords{for all $r\in R,\ a\in A,\ b\in B$}
\end{equation}
This is equivalent to saying that homomorphism 
$\lambda_\ast: A \rightarrow \Hom(B,C)$ induced by $\lambda$ is an $R$-module homomorphism (where $R$ acts on  $\Hom(A,C)$ by pre-composition with its action on $A$) or indeed to the same thing with $A$ and $B$ interchanged. 
We shall summarise this situation by saying that $\cH=\HeiStan$ \textit{`has $R$-structure'}. The phrase should strictly be applied to the tuple $(A,B,C,\lambda)$ defining $\cH$, moreover the $R$-structure is usually \textit{not} determined by the group structure of $\cH$ itself. Nevertheless, in this situation, identifying $\cH/m(C)$ with $A\oplus B$ imposes on it the natural structure of an $R$-module. 

Suppose that   
$\t \cH=\Hei{\t A}{\t B}{\t C}{\t \lambda}$ also has $R$-structure and  
$\theta:\cH\rightarrow\t \cH$ is any homomorphism of groups with 
$\theta(m(C))\subset m(\t C)$. (For instance, if $\theta$ is surjective and  condition ND-C is satisfied by both tuples, so $m(C)=Z(\cH)$ and $m(\t C)=Z(\t \cH)$.) 
Then we shall say abusively that $\theta$ is $R$-linear (or an $R$-homomorphism) iff the induced homomorphism 
$\bar \theta  :A\oplus B\rightarrow \t A\oplus \t B$ is $R$-linear in the usual sense. In particular,  
a diagonal homomorphism $t_A\times t_B\times t_C$ is $R$ linear iff $t_A$ and $t_B$ are $R$-linear. If condition ND-C is satisfied, we may  define $\Aut_R(\cH)$ to be the
subgroup of \textit{$R$-linear} automorphisms in $\Aut(\cH)$ and $\Aut^0_R(\cH)$ to be $\Aut_R(\cH)\cap\Aut^0(\cH)$ which contains $\Inn(\cH)$ since the latter acts by the identity on $A\oplus B$. Similarly, we set $\Out_R(\cH):=\Aut_R(\cH)/\Inn(\cH)<\Out(\cH)$ and $\Out_R^0(\cH):=\Aut_R^0(\cH)/\Inn(\cH)<\Out^0(\cH)$. It seems
 natural to investigate $\cA$-regular bouquets for $\cH$ with  $\cA$ contained in $\Out_R(\cH)$ (and with respect to $\sigma_p$, assuming the latter is irreducible).
\subsection{GHGs of Direct Sums and Tensor-Product Representations}
Consider now the situation where, for given (cyclic) $C$,  the tuples $(A_i,B_i,C,\lambda_i)$ for $i=1,\ldots,m$  all satisfy condition ND-C, with $A_i$, and $B_i$ being $R$-modules and $\lambda_i$ being $R$-balanced for all $i$. Suppose $A$ is the direct sum $\bigoplus_{i=1}^m A_i$ as an $R$-module and similarly 
$B=\bigoplus_{i=1}^m B_i$ as an $R$-module. We define a pairing $\lambda$: 
$$
\displaynomapdef
{\lambda:A\times B}
{C}
{((a_1,\ldots,a_m),(b_1,\ldots,b_m))}
{\lambda_1(a_1,b_1)+\ldots+\lambda_1(a_m,b_m)}
$$ 
It is easy to see that $\lambda$ is $R$-balanced, bilinear and non-degenerate, so $(A,B,C,\lambda)$ satisfies~condition ND-C and $\HeiStan$ is a generalised Heisenberg group with $R$-structure. 
For each $i$, let 
$j_{A_i}$ 
and $j_{B_i}$ be the natural
injective $R$-module maps from $A_i$ into $A$ and
$B_i$ into $B$, respectively. By 
Proposition~\ref{prop: criterion for diagonal homs} we have an (injective, $R$-linear) diagonal homomorphism for each $i$
$$
\ut_i:=j_{A,i}\times j_{B,i}\times \id_C : \cH_i:=\cH(A_i,B_i,C,\lambda)
\longrightarrow \cH:=\HeiStan
$$
For each $(i,k)$ with $i\neq k$ we have $\lambda(j_{A,i}(A_i),j_{B,k}(B_k))=\{0\}$. Thus  $\im(\ut_i)\cong\cH_i$ and $\im(\ut_k)\cong\cH_k$ centralise each other and, moreover,   
$\im(\ut_i)\cap\im(\ut_k)=m(C)=Z(\cH)$. One proves easily:
\begin{prop}\label{prop:direct product of Hei}{\bf (Heisenberg Group of Direct Sums under Condition ND-C)}\ \\
	The map from the direct product 
	$$
	\displaymapdef{\theta}{\prod_{i=1}^m\cH_i}{\cH}{(h_1,\ldots,h_m)}{\ut_1(h_1)\ldots\ut_m(h_m)}
	$$
	is a surjective group homomorphism and $\im(\ut_i)$is normal in $\cH$  for all $i$. Moreover 
$$
\ker\theta=\{(m(c_1),\ldots,m(c_m))\,:\, c_1+\ldots+c_m=0\}<
\prod_{i=1}^m Z(\cH_i)=Z\left(\prod_{i=1}^m\cH_i\right)
=\theta\inv(Z(\cH))
$$ 
%
In particular, $\theta$ induces an isomorphism modulo centres.\ePf  
\end{prop}
Now suppose $p:C\rightarrow \bbC^\times$ is an injective homomorphism so that  $\sigma_p:\cH\rightarrow\U(\cM(A))$ is  faithful and SV. It follows (\eg\ by Proposition~\ref{prop:direct product of Hei} and Lemma~\ref{lemma:when rho and tilderho are SV}) 
that the inflation of $\sigma_p$ to the group $\prod_{i=1}^m\cH_i$ -- which is class $\leq 2$ nilpotent, but does not have cyclic centre if $m>1$ -- is also SV. Explicitly, this is the representation $\rho_1:=\sigma_p\circ\theta$ so that
$$
  \underline{h}\cdot_{\rho_1}f=\theta(\underline{h})\cdot_{\sigma_p}f\addwords{\ \ \  $\forall \underline{h}\in\prod_{i=1}^m\cH_i$,\ $\forall f\in \cM(A)$} 
$$
On the other hand, writing $\sigma_{i,p}$ for the Schr\"odinger representation of $\cH_i$ on $\cM(A_i)$, we have the tensor-product representation of $\prod_{i=1}^m\cH_i$ on 
$\bigotimes_{i=1}^m\cM(A_i)$ namely the representation
$$
\rho_2:=\bigotimes_{i=1}^m\sigma_{i,p}
\ \ \ \mbox{so}\ \ \ 
\underline{h}\cdot_{\rho_2}
\left(
\sum_{s=1}^S f_{1,s}\otimes\ldots\otimes f_{m,s}
\right)
=
\sum_{s=1}^S ( h_1\cdot_{\sigma_{1,p}} f_{1,s})\otimes\ldots\otimes (h_m\cdot_{\sigma_{m,p}} f_{m,s}) 
$$
where $\underline{h}=(h_1,\ldots,h_m)\in\prod_{i=1}^m\cH_i$ and $f_{i,s}\in \cM(A_i)$ for all $i$ and $s$. (All tensor products are, of course, over $\bbC$.) Note also that $\rho_2$ is unitary w.r.t.\ the PDHF on 
$\bigotimes_{i=1}^m\cM(A_i)$ given by 
$$
\langle
f_1\otimes\ldots\otimes f_m, 
g_1\otimes\ldots\otimes g_m 
\rangle:=\langle f_1,g_1 \rangle\ldots\langle f_m,g_m \rangle\in \bbC
$$
(extended bilinearly) and that we have a linear isomorphism
$$
\displaymapdef{\cV}{\bigotimes_{i=1}^m\cM(A_i)}{\cM(A)}
{f_1\otimes\ldots\otimes f_m}
{\left (f:(a_1,\ldots,a_m)\mapsto f_1(a_1)\ldots f_m(a_m)\right)}
$$ 
which sends orthonormal basis elements  $e_{a_1}\otimes\ldots\otimes  e_{a_m}$ to orthonormal basis elements $e_{(a_1,\ldots,a_m)}$ and so respects the PDHFs. 
\begin{prop} In the above situation and notations, the map $\cV$ defines an isomorphism between 
$\rho_2$ and $\rho_1$ as unitary representations of 
$\prod_{i=1}^m\cH_i$, that is:
\begin{equation}\label{eq:cV a unitary iso from rho2 to rho1}
\cV(\underline{h}\cdot_{\rho_2} v)
=
\underline{h}\cdot_{\rho_1}\cV(v)
\addwords{ $\forall \underline{h}\in
	\prod_{i=1}^m\cH_i$,  
	$\forall v=\sum_{s=1}^S f_{1,s}
	\otimes\ldots\otimes f_{m,s}\in \bigotimes_{i=1}^m\cM(A_i)$}
\end{equation}
In particular, $\rho_2$ is an SV-representation and factors through a representation $\bar{\rho}_p$, say, of $\cH(A,B,C;\lambda)$ on 
$\bigotimes_{i=1}^m\cM(A_i)$ which is isomorphic by $\cV$ to $\sigma_p$, as unitary representations of $\cH$. Explicitly, we can write  
\begin{eqnarray*}
h((a_1,\ldots,a_m),(b_1,\ldots,b_m), c)\cdot_{\bar{\rho}_p}
	\left(
	\sum_{s=1}^S f_{1,s}\otimes\ldots\otimes f_{m,s}
	\right)
&=&p(c)\sum_{s=1}^S ( h'_1\cdot_{\sigma_{1,p}} f_{1,s})\otimes\ldots\otimes (h'_m\cdot_{\sigma_{m,p}} f_{m,s})\\
\end{eqnarray*}
where $h'_i=h(a_i,b_i,0)\in \cH_i$ for $i=1,\ldots,m$.
\end{prop}
\bPf\ The main thing to check is~(\ref{eq:cV a unitary iso from rho2 to rho1}) and, for this, it suffices to take  $\underline{h}$ to be of form $(1,\ldots,1,h(a,b,c),1,\ldots,1)$ with $a\in A_i,\ b\in B_i,\ c\in C$ for some $i\in\{1,\ldots,m\}$ 
and $v$ to be a basis element 
of form $e_{a_1}\otimes\ldots\otimes  e_{a_m}$ where $a_i\in A_i\ \forall i$. In this case, we find that both sides of~(\ref{eq:cV a unitary iso from rho2 to rho1}) are equal to 
$p(\lambda_i(a_i-a,b)+c)e_{(a_1,\ldots,a_i-a,\ldots,a_m)}$. The rest is LTR.
\ePf
\subsection{Diagonal $R$-Isomorphisms under Condition ND-C}
We suppose now that $(A,B,C,\lambda)$ and 
$({\t A},{\t B},{\t C},{\t \lambda})$ are two tuples, both once again satisfying condition ND-C and such that 
the groups $\HeiStan$  and $\cH({\t A},{\t B},{\t C},{\t \lambda})$ have $R$-structure for the same commutative ring $R$ (possibly $R=\bbZ$). In this situation, the following result shows that diagonal \textit{$R$-isomorphisms} of the form $t_A \times t_B \times t_C$ exist and are determined uniquely by the pairs $(t_A,t_C)$, with $t_A\in\Hom_R(A,\t A)$ and $t_C\in\Hom(C,\t C)$. (The same is not true without the isomorphism condition.) This has important consequences for the isomorphism of Heisenberg groups and  
their Schr\"odinger representations. (An analogous result  holds for pairs $(t_B,t_C)$ and can proven by an analogous argument.)
\begin{thm}{\bf (Extension to Diagonal Isomorphisms)}\label{prop:Extension to Diagonal Isomorphisms}\\
	Suppose that $(A,B,C,\lambda)$, 
	$({\t A},{\t B},{\t C},{\t \lambda})$  are as specified above and suppose that $t_A\in\Hom_R(A,\t A)$ and $t_C\in\Hom(C,\t C)$ are both \emph{isomorphisms}. 
	\begin{enumerate}
		\item\label{part:Extension to Diagonal Isomorphisms1} There exists a unique element 
		$t_B$ of $\Hom(B,\t B)$ such that $\underline{t}:=t_A\times t_B\times t_C$ is a homomorphism from $\HeiStan$  to $\cH({\t A},{\t B},{\t C},{\t \lambda})$.
		\item\label{part:Extension to Diagonal Isomorphisms2} Moreover, $t_B$ and $\underline{t}$ are $R$-isomorphisms.
	\end{enumerate} 
\end{thm}
\bPf\ We have $R$-homomorphisms 
$
\lambda_\ast: B \rightarrow \Hom(A,C)
$ and $
\t \lambda_\ast: \t B \rightarrow \Hom(\t A,\t C)
$
induced by $\lambda$ and $\t \lambda$ respectively. 
By Proposition~\ref{prop: criterion for diagonal homs}, given $t_B\in\Hom(B,\t B)$, the map $\underline{t}$ will be a homomorphism from $\HeiStan$  to $\cH({\t A},{\t B},{\t C},{\t \lambda})$ iff 
$
t_C\circ\lambda_\ast(b)=
{\t \lambda}_\ast(t_B(b))\circ 
t_A
$
in $\Hom(A,\t C)$, for all $b\in B$. But $t_A$ is an isomorphism and we have seen that condition ND-C for 
$({\t A},{\t B},{\t C},{\t \lambda})$ implies that ${\t \lambda}_\ast$ is also an isomorphism so this condition can be rewritten as 
\begin{equation}\label{eq: criterion for t_B}
{\t \lambda}_\ast\inv\left(t_C\circ\lambda_\ast(b)\circ t_A\inv\right)=
t_B(b)\end{equation}
for all $b\in B$. (The expression in large parentheses on the L.H.S.\ is clearly an element of $\Hom(\t A,\t C)$.) Thus part~\ref{part:Extension to Diagonal Isomorphisms1} will follow if we can show that the L.H.S.\ of~(\ref{eq: criterion for t_B}), as a function of $b\in B$, is a homomorphism from $B$ to $\t B$. But this follows from the fact that $\lambda_\ast$, $t_C$ and ${\t \lambda}_\ast\inv$ are homomorphisms. For part~\ref{part:Extension to Diagonal Isomorphisms2}, it suffices to show that $t_B$, as defined by~(\ref{eq: criterion for t_B}), is an $R$-isomorphism, as this will imply the same for $\underline{t}$. The $R$-linearity of $t_B$ follows from that of 
 $\lambda_\ast$, $t_A\inv$ and ${\t \lambda}_\ast\inv$ so it only remains to show that $t_B$ is bijective. Since $t_A$ is an isomorphism we have $|A|=|\t A|$, \ie\  $|B|=|\t B|$ by condition ND-C, so injectivity suffices. But for $b,b'\in B$, the equality
 $$
 {\t \lambda}_\ast\inv\left(t_C\circ\lambda_\ast(b)\circ t_A\inv\right)=
 {\t \lambda}_\ast\inv\left(t_C\circ\lambda_\ast(b')\circ t_A\inv\right)
 $$
implies $b=b'$ since $t_C$ is injective and $\lambda_\ast$ is too, by condition ND-C for 
$(A,B,C,\lambda)$.\ePf
\begin{cor}\label{cor:implication about iso}
	Suppose $(A,B,C,\lambda)$ and 
	$({\t A},{\t B},{\t C},{\t \lambda})$ satisfy condition ND-C and 
	the groups $\HeiStan$  and $\cH({\t A},{\t B},{\t C},{\t \lambda})$ have $R$-structure. Then
	\begin{equation}\label{eq:implication about iso}
		\mbox{
			($A \cong\t A$ 
			as $R$-modules and $C \cong \t C$) $\Longrightarrow$
			($\HeiStan\cong\cH({\t A},{\t B},{\t C},{\t \lambda})$ over $R$)
		}\ \ \Box\ \ \	\end{equation}
\end{cor}
\rem\  Note in particular that, if condition ND-C is satisfied, then the $R$-isomorphism type of $\HeiStan$ doesn't depend on $\lambda$, only on that of $A$ (as well the isomorphism type of $C$, but that is determined by $r=|C|$)! Obviously, implication~(\ref{eq:implication about iso}) holds equally with $B\cong\t B$ in place of $A\cong\t A$ on the L.H.S., by a similar argument. Also the converse implication also holds under condition ND-C, at least if $R=\bbZ$: if $\HeiStan$ and  $\cH({\t A},{\t B},{\t C},{\t \lambda})$ are isomorphic as groups, then so are their centres and quotients-by-centres. Thus $C\cong\t C$ and $A\oplus B\cong \t A\oplus \t B$, hence   
$A\oplus A\cong \t A\oplus \t A$ so $A\cong\t A$ by the structure theorem for finite abelian groups.\bigskip\\
Let us take $(A,B,C,\lambda)=(\t A,\t B,\t C,\t \lambda)$ and $t_C=\id_C$ in Theorem~\ref{prop:Extension to Diagonal Isomorphisms}. Writing $\alpha$ for $t_A\in \Aut_R(A)$, we easily deduce: 
\begin{cor}{\bf (Diagonal Automorphisms)}\label{cor:Diagonal Automorphisms}\\
	Suppose that $(A,B,C,\lambda)$ satisfies condition ND-C and $\cH:=\HeiStan$ has $R$-structure. Then there is an injective homomorphism
	$$
	\displaymapdef{\Delta}{\Aut_R(A)}{\Aut^0_R(\cH)}{\alpha}{\alpha \times \beta \times \id_C}
	$$
	where $\beta$ is the unique element of $\Aut(B)$ satisfying
	\begin{equation}\label{eq: relation between alpha and beta in diag aut}
	\lambda(\alpha(a),\beta(b))=\lambda(a,b)\ \forall a,b. 
	\end{equation}
	 In fact, $\beta$ lies in $\Aut_R(B)$ and is the inverse of the `dual  automorphism' to $\alpha$ w.r.t.\ $\lambda$. In particular, if $\alpha$ is the action of $r\in R ^\times$ then $\beta$ is the action of $r\inv$.)\ePf  
\end{cor}
The image $\Delta(\Aut_R(A))$ of $\Delta$ in $\Aut^0_R(\cH)$ will be called the subgroup of \textit{diagonal ($R$-) automorphisms of $\cH$ (fixing the centre)}. 
Combining Theorem~\ref{prop:Extension to Diagonal Isomorphisms} and Proposition~\ref{prop:"for future use"} gives:
\begin{cor}\label{cor:Diag Auts}
Suppose that 
$(A,B,C,\lambda)$, 
$({\t A},{\t B},{\t C},{\t \lambda})$, $t_A$ and $t_C$ satisfy the conditions of Theorem~\ref{prop:Extension to Diagonal Isomorphisms} and that $\underline{t}=t_A\times t_B\times t_C$ is the unique $R$-isomorphism described therein. Suppose also  $p:C\rightarrow\bbC^\times$ and ${\t p}:{\t C}\rightarrow\bbC^\times$ are homomorphisms satisfying ${\t p}\circ t_C=p$ and let $\sigma_p:\cH:=\HeiStan\rightarrow\cM(A)$ and 
$\sigma_{\t p}:\t \cH:=\Hei{\t A}{\t B}{\t C}{\t \lambda}\rightarrow\cM(\t A)$ be the (left) Schr\"odinger representations associated to $p$ and $\t p$ respectively. Then 
$$
\displaymapdef{\check t_A}{\cM(A)}{\cM(\t A)}{f}{f\circ t_A\inv}
$$ satisfies equation~(\ref{eq:actions commute with projection})
and is an isomorphism of unitary representations from $\sigma_p$ to $\sigma_{\t p}$. 
\end{cor}
\bPf\ It only remains to note that $\check t_A$ is an isomorphism (obvious) and that it preserves PDHFs (by Remark~\ref{rem:injectivity of tAstar}).\ePf
\noindent Corollary~\ref{cor:Diag Auts} now yields the Weil representation on $\Delta(\Aut(A))<\Aut^0(\HeiStan)$: 
\begin{cor}\label{cor: Weil Rep for Diag Auts}
	{\bf (Weil Representation for Diagonal Automorphisms in $\Aut^0(\cH)$)}\\
Suppose that $(A,B,C,\lambda)$ satisfies condition ND-C and  that $p:C\rightarrow\bbC^\times$ is a homomorphism such that the Schr\"odinger representation $\sigma_p$ of $\cH=\HeiStan$ is irreducible (\eg\ $p$ is injective). Thus the Weil representation $\cT_{\sigma_p}:\Aut^0(\cH)\rightarrow {\rm PU}(\cM(A))$ with respect to $\sigma_p$ is well-defined. Then, for any $\alpha\in \Aut(A)$ the element 
$
\cT_{\sigma_p}(\Delta(\alpha))
$ 
of ${\rm PU}(\cM(A))$ is precisely the image of $\check \alpha\in \U(\cM(A))$, namely the map sending $f\in\cM(A)$  to $\check \alpha(f)=f\circ\alpha\inv$. 
\end{cor}
\bPf\ Take $R=\bbZ$, 
$({\t A},{\t B},{\t C},{\t \lambda})=(A,B,C,\lambda)$,  $t_A=\alpha$ and  $t_C=\id_C$ in Corollary~\ref{cor:Diag Auts}, so that $\t p=p$. Then $\underline{t}$ equals $\Delta(\alpha)$ by definition of the latter and $\check t_A=\check\alpha$. So, according to the last Corollary, $\check \alpha$ lies in $\U(\cM(A))$ and the following version of equation~(\ref{eq:actions commute with projection}) holds 
\begin{equation}\label{eq:first of pf of  Weil Rep for Diag Auts}
	\Delta(\alpha)(h)\cdot_{\sigma_{p}} \check \alpha(f) =  \check \alpha (h\cdot_{\sigma_p} f)
\end{equation}
for all $h\in\HeiStan$ and $f\in\cM(A)$. In other words 
$$
 \check \alpha\circ\sigma_p(h)=\sigma_p(\Delta(\alpha)(h))
\circ \check \alpha
$$ 
in $\U(\cM(A))$ for all $h\in\HeiStan$.
Comparing this last equation with~(\ref{eq:def T_phi})  (with $\rho=\sigma_p$ and $g=h\in\cH$) and the definition~(\ref{eq:def of Weil Rep}) of $\cT_{\sigma_p}$  
gives the result.\ePf  
\noindent By contrast to Corollary~\ref{cor: Weil Rep for Diag Auts} the explicit characterisation of $\cT_{\sigma_p}$ on  \textit{non}-diagonal elements of 
$\Aut^0(\cH)$ is a much harder -- but also more interesting -- task, even in the Base Case. See Remark~\ref{rem:TheCliffordGroupandTheExtCG}.\bigskip\\
\rem\ A special case occurs when $\HeiStan$ has $R$-structure and $\alpha$ is the action of $r\in R^\times$. In particular, if we take  $r=-1$, so that $\alpha=\alpha\inv$ is inversion in $A$, we find  $\Delta(\alpha)=\phi_{-1}$. Thus (if $\sigma_p$ is irreducible) the  
Corollary 
tells us that $\cT_{\sigma_p}(\phi_{-1})$ sends $f(x)\in{\cM(A)}$ to $f(-x)$.\bigskip\\
\rem\label{rem:R-lin extn of Weil}\ {\bf (An $\bbR$-Linear Extension of the Weil Representation)}
Suppose $r \neq 2$ so that the automorphism $\phi^-$ of $\cH$ sending $h(a,b,c)$ to $h(a,-b,-c)$ is of order exactly $2$ and does \textit{not} lie in $\Aut^0(\cH)$. We define $\bc\in \GL_\bbR(\cM(A))$ of order two by $\bc(f)(a)=f(a)^\ast$ (complex conjugate) for all $f\in\cM(A)$ and $a\in A$. One checks easily the following analogue of~(\ref{eq:first of pf of  Weil Rep for Diag Auts}):
$$
	\phi^-(h)\cdot_{\sigma_{p}} \bc(f) =  \bc  (h\cdot_{\sigma_p} f)
$$
for all $h$ and $f$. Now, given any $\phi\in\Aut^0(\cH)$ and $T_\phi\in\U(\cM(A))$ representing $\cT_{\sigma_p}(\phi)$ modulo $\bbT$, it follows easily  from~(\ref{eq:def T_phi}) (with $\rho=\sigma_p$ and $g=h\in\cH$) and the last equation (twice) that   
$$
\phi^-\phi\phi^-(h)\cdot_{\sigma_{p}} \bc T_\phi \bc(f) =  \bc T_\phi \bc  (h\cdot_{\sigma_p} f)
$$
for all $h$ and $f$. But $\phi^-\phi\phi^-$ and $\bc T_\phi \bc$ clearly lie in $\Aut^0(\cH)$ and $\U(\cM(A))$ respectively, so, by the definition of $\cT_{\sigma_p}$, we must have 
\begin{equation}\label{eq:to extend cT}
	\cT_{\sigma_p}(\phi^-\phi\phi^-)=\bc \cT_{\sigma_p}(\phi)\bc
	\pmod{\bbT}
\end{equation}
in  ${\rm PU}(\cM(A))$. Let us denote by $\widetilde{\Aut^0}(\cH)$ the subgroup of $\Aut(\cH)$ consisting of those automorphisms acting on $Z(\cH)\cong\Zover{r}$ by $\pm 1$. Then, clearly, $\widetilde{\Aut^0}(\cH)=\Aut^0(\cH)\stackrel{\cdot}{\cup} \phi^-\Aut^0(\cH)$ and it follows easily from~(\ref{eq:to extend cT})  that $\cT_{\sigma_p}$ extends to a (unique)  \textit{homomorphism} $\t \cT_{\sigma_p}:\widetilde{\Aut^0}(\cH)\rightarrow\GL_\bbR(\cM(A))/\bbT$ sending $\phi^-$ to  $\bc$ ($\bmod\ \bbT$). By construction, it clearly satisfies 
$$
	\t \cT_{\sigma_p}(\upsilon)\circ\sigma_p(h)= \sigma_p(\upsilon(h))\circ \t \cT_{\sigma_p}(\upsilon)
	\pmod{\bbT}
$$
for all $\upsilon\in\widetilde{\Aut^0}(\cH)$ and $h\in \cH$. Note that $\phi_-$ lies in $\widetilde{\Aut^0}(\cH)$ and equals $\phi^-\phi_{-1}$. Thus, by the previous remark, $\t \cT_{\sigma_p}(\phi_-)$  is the image in $\GL_\bbR(\cM(A))/\bbT$  of the map sending $f$ to $\bc(f)$ composed with inversion in $A$. We return to the  representation
$\t \cT_{\sigma_p}$ in the Base Case in \S\ref{subsec:the Base Case}.
\section{Automorphism Groups and Displacement Operators}\label{sec:Aut Gps}
\subsection{ Hypotheses, Notation and Orientation} 
Throughout this section, we assume that $\cH:=\HeiStan$ is a generalised Heisenberg group with $R$-structure for a commutative ring $R$ (possibly $R=\bbZ$), as defined in the previous section. In particular $A$ and $B$ are $R$-modules, so also is $\AdsB$. We shall suppose that condition ND-C holds for the tuple $(A,B,C,\lambda)$. Thus  $Z\cong C$ is cyclic of order $r$ and $A$ and $B$ are dual via $\lambda$.
We shall further assume for simplicity 
that $\cH$ is non-abelian, \ie\ $A$ is non-trivial, hence also $C$. Thus $r>1$ and $\cH\neq Z\neq \{1\}$. We shall usually write $\ba$ generically for $(a,b)$,  $\ba'$ for $(a',b')$ \textit{etc} in $A\oplus B$. We shall often  write $h(\ba,c)$ for $h(a,b,c)$ \textit{etc} in $\cH$. To further save notation, condition ND-C allows us to identify $\bar{\cH}:=\cH/Z$ with $A\oplus B$ by means of $\pi$ (\textit{i.e.} we identify $\overline{h(\ba,c)}$ with $\ba$ for any $c\in C$) when this is convenient. In particular,  the non-degenerate,  skew-symmetric bilinear form $\delta_\cH :\bar{\cH}\times \bar{\cH}\rightarrow Z$ (\textit{cf}~(\ref{eq:def of deltaG for any G class 2 nilp})) then becomes $m\circ \delta$ where 
\begin{equation}\label{eq:def of delta lambda}
	\displaymapdef{\delta=\delta_\lambda}{(A\oplus B)\times (A\oplus B)}{C}{(\ba,\ba')}{\lambda(a,b')-\lambda(a',b)}
\end{equation}
(\textit{cf}~(\ref{eq:def of deltaG for any G class 2 nilp})  and (\ref{eq: commutators in HeiStan})) and since $\lambda$ is $R$-balanced, so also is  $\delta$.

Two longer term goals are, firstly, to render explicit the Weil representation of $\Aut^0_R(\cH)$ w.r.t.\ an irreducible Schr\"odinger representation $\sigma_p$ of $\cH$, and, secondly, to investigate invariant and regular $\cH$-bouquets w.r.t.\ $\sigma_p$ and subgroups $\hat \cE$ and $\cA$ of $\Out^0_R(\cH)$ respectively. An essential first step for both of these is to get a description of $\Aut^0_R(\cH)$and $\Out^0_R(\cH)$ in terms of $A$, $B$ and $\lambda$. We achieve this in the first subsection here under the hypothesis that \textit{$r$ is odd}. This  allows us to define a map $\cD:\cH/Z=A\oplus B \rightarrow\cH$ which is right inverse to $\pi$ and has special properties:
its associated $2$-cocycle is $\frac12 \delta_\cH$ (see~(\ref{eq:D is a hom mod Z})) and it `commutes' with powers and with diagonal automorphisms  (see~(\ref{eq:D commutes w powers}) and~(\ref{eq:D commutes w diag autoorphisms})). 
The composition of $\cD$ with the Schr\"odinger representation yields elements $D(\ba)$ of $\U(\cM(A))$ which form a $\bbC$-basis of $\End_\bbC(\cM(A))$. In the Base Case these elements are the unitary endomorphisms of $\cM(\Zover{d})\cong\bbC^d$ referred to as \textit{`displacement operators'} in most of the SIC-related literature of quantum physics, from which the name, of course, derives (see \S\ref{subsubsection:base case odd d}). 
In the second subsection here we briefly explain the connection between the overlap-maps introduced in \S\ref{sec:Overlap- and Angle-Maps, and Bouquets} -- now specialised by taking $\cR=\cD$ -- 
and the coefficients of the $D(\ba)$'s in certain projectors in $\End_\bbC(\cM(A))$. (See Remark~\ref{rem:  angle-maps related to projectors}.) In the Base Case, these coefficients are (essentially) the so-called \textit{`overlaps'} -- another word with a  `physical' etymology giving rise to the terminology used in this paper. 

\subsection{Automorphism Groups \textit{via} the Map $\cD$}
\label{subsec:Aut Gps via D}
To start, we do not suppose $r$ to be odd. Recall from Subsection~\ref{subsec:Reg and Inv} the homomorphism  
$$
\displaynomapdef{\Aut(\cH)}{\Aut(\bar \cH)=\Aut( \AdsB)}{\phi}{\bar\phi}
$$
factoring through the quotient map $\Aut(\cH)\rightarrow\Out(\cH)$ which sends $\phi$ to $\hat\phi$. Now suppose that $\phi$ lies in $\Aut^0_R(\cH)$ \ie\ it is $R$-linear (modulo $Z$) and fixes $Z$ elementwise. Then, for any $h,h'\in\cH$ we have:
$$
\delta_\cH(\bar h, \bar h')=
\phi\left(\delta_\cH(\bar h, \bar h')\right)=
\phi([h, h'])=
[\phi(h), \phi(h')]=
\delta_\cH(\bar\phi(\bar h), \bar\phi(\bar h'))
$$ 
Writing $\delta_H$ as $m\circ\delta$ and $\ba,\ba'$ for ${\bar h}, {\bar h}'$ respectively, we have shown that $\bar\phi$ lies in the \textit{symplectic group} $\SpABd{R}$ w.r.t.\ $\delta$, that is:
\begin{equation}\label{eq:defn SpABd}
	\SpABd{R}:=\{\bar\phi\in\Aut_R(\AdsB)\,:\, 
	\delta(\ba,\ba')=\delta(\bar{\phi}(\ba),\bar{\phi}(\ba'))\ \mbox{for all $\ba,\ba'\in\AdsB$}\}
\end{equation}
and 
\begin{prop}\label{prop:the hom bar Theta} There is a well defined homomorphism 
$$
\displaymapdef{\bar\Theta}{\Out^0_R(\cH)}{\SpABd{R}}{\hat \phi}{\bar \phi \ \ \ \ \ \Box}
$$	
\end{prop}
We now make the 
\begin{hyp}\label{hyp:r is odd}
	\mbox{$r$ is odd}
\end{hyp}
(This implies that  $e=\exp(\cH/Z)$ and  $s=|A|=\sqrt{|\cH/Z|}$ are also odd, so it is equivalent to $2\ndiv |\cH|$.) Assuming this hypothesis, Theorem~\ref{thm: 3x2 diagram of aut gps} will show that $\bar\Theta$ is an isomorphism and, moreover, give a full characterisation of $\Aut^0_R(\cH)$. 

To start, Hypothesis~\ref{hyp:r is odd} implies that every $c\in C$ has a unique square-root, written $\frac12 c$ so we can define a map $\cD:A\oplus B \rightarrow\cH$ by
\begin{equation}\label{def:displacement operator}
	\cD(\ba):=h\left(\ba,\frac12\lambda(a,b)\right)
\end{equation}
Thus $\pi\circ \cD=\id_{\AdsB}$ and $\cD(0)$ is the identity. Direct calculation gives the important relation 
\begin{equation}\label{eq:D is a hom mod Z}
	\cD(\ba)\cD(\ba')=
	\cD(\ba+\ba')m\left(\half\delta(\ba,\ba')\right)
\end{equation}
by~(\ref{eq:def of delta lambda}), and the skew-symmetry of 
$\delta$ then yields the special property
\begin{equation}\label{eq:D commutes w powers}
\cD(\ba)^n = \cD(n\ba)\addwords{for all $n\in\bbZ$}
\end{equation}
(It follows from this that the exponent of $\cH$ is exactly $r$ under Hypothesis~\ref{hyp:r is odd}.) 
For each $\phi\in\Aut^0_R(\cH)$, we define $\eta_\phi:\AdsB\rightarrow C$ by 
\begin{equation}\label{eq:phi(D(a))}
	\phi(\cD(\ba))=\cD(\bar{\phi}(\ba))m(\eta_\phi(\ba))\addwords{for all $\ba\in\AdsB$}
\end{equation}
In fact, $\eta_\phi$ must be a homomorphism of abelian groups. To see this,  we apply $\phi$ to equation~(\ref{eq:D is a hom mod Z}) with the two sides interchanged and use~(\ref{eq:phi(D(a))}) on each side to get 
\begin{eqnarray*}
	\cD(\bar{\phi}(\ba+\ba'))m\left(\eta_\phi(\ba+\ba')+\half\delta(\ba,\ba')\right)&=&\cD(\bar{\phi}(\ba))m(\eta_\phi(\ba))
	\cD(\bar{\phi}(\ba'))m(\eta_\phi(\ba'))\\
	&=&\cD(\bar{\phi}(\ba))
	\cD(\bar{\phi}(\ba'))m\left(\eta_\phi(\ba)+\eta_\phi(\ba')\right)\\
	&=&\cD(\bar{\phi}(\ba)+\bar{\phi}(\ba'))
	m\left(\frac12 \delta(\bar{\phi}(\ba),\bar{\phi}(\ba'))+\eta_\phi(\ba)+\eta_\phi(\ba')\right)\\	
	&=&\cD(\bar{\phi}(\ba+\ba'))
	m\left(\frac12 \delta(\ba,\ba')+\eta_\phi(\ba)+\eta_\phi(\ba')\right)\\
\end{eqnarray*}
using~(\ref{eq:D is a hom mod Z}) again and the fact that  $\bar{\phi}$ lies in $\SpABd{R}$. Cancelling and using the injectivity of $m$, gives 
$\eta_\phi(\ba+\ba')=\eta_\phi(\ba)+\eta_\phi(\ba')$, as required, so we now have a well-defined map of \textit{sets}
$$
\displaymapdef{\t \Theta_\cD}{\Aut^0_R(\cH)}{\Hom(\AdsB,C)\times\SpABd{R}}{\phi}{(\eta_\phi,\bar{\phi})}
$$
If ${\t \Theta_\cD}(\phi)={\t \Theta_\cD}(\phi')$ then, by~(\ref{eq:phi(D(a))}), $\phi$ and $\phi'$ agree on $\cD(\ba)$ for all $\ba\in\AdsB$ and since they both fix $Z$, they are equal, so ${\t \Theta}$ is injective. To see that ${\t \Theta}$ is surjective, we suppose given  $\eta\in\Hom(\AdsB,C)$ and $\alpha\in\SpABd{R}$. Each element of $\cH$ can be written uniquely as $\cD(\ba)m(c)$ with $\ba\in\AdsB$ and $c\in C$, so we can define a map 
\begin{equation}\label{eq:defn of nu}
\displaymapdef{\nu}{\cH}{\cH}{\cD(\ba)m(c)}{\cD(\alpha(\ba))m(\eta(\ba)+c)}	
\end{equation}
It suffices to show that $\nu$ lies in $\Aut^0_R(\cH)$. (For then, taking $c=0$ and comparing with~(\ref{eq:phi(D(a))})
shows that $\eta_\nu=\eta$ and $\bar{\nu}=\alpha$ so  ${\t \Theta_\cD}(\nu)=(\eta,\alpha)$.) 
Now, the injectivity of $\cD$, $m$ and $\alpha$ implies that of $\nu$, so $\nu$ is bijective since $\cH$ is 
finite. Also, taking $\ba=0$ in~(\ref{eq:defn of nu}) shows that $\nu$ fixes $Z$. To prove that $\nu$ is a homomorphism, we first check that for any $\ba,\ba'\in\AdsB$ we have, by~(\ref{eq:defn of nu}) and~(\ref{eq:D is a hom mod Z}), the two equations 
\begin{eqnarray*}
\nu(\cD(\ba))\nu(\cD(\ba'))&=&\cD(\alpha(\ba)+\alpha(\ba'))m\left(\frac12\delta(\alpha(\ba),\alpha(\ba'))+\eta(\ba)+\eta(\ba')\right)\\
\nu(\cD(\ba)\cD(\ba'))&=&\cD(\alpha(\ba+\ba'))m\left(\eta(\ba+\ba')+\half\delta(\ba,\ba')\right)
\end{eqnarray*}
Their RHS's are equal since $\alpha$ and $\eta$ are homomorphisms and $\alpha$ preserves $\delta$. Since also $\nu(hm(c'))=\nu(h)m(c')$ for all $h\in\cH,c'\in C$ by~(\ref{eq:defn of nu}), it follows that $\nu$ is a homomorphism. Thus $\nu$ lies in $\Aut^0_R(\cH)$. Thus ${\t \Theta_\cD}$ is surjective, hence a bijective map of sets. 

One way to turn $\t \Theta_\cD$ into a group isomorphism is to simply carry the group structure across from $\Aut^0_R(\cH)$
to the product set $\Hom(\AdsB,C)\times\SpABd{R}$. However, this results in an unconventional 
group structure on the latter. Here is one way to tackle this (essentially notational) problem. A more conventional group structure on the product set is obtained  by first making  $\SpABd{R}$ act 
\textit{on the left} on 
$\Hom(\AdsB,C)$, namely 
$\bar \phi\in\SpABd{R}$ sends 
$\eta\in \Hom(\AdsB,C)$ to 
$\eta\circ\bar \phi\inv$. 
The usual definition of the 
semidirect product with this action (see e.g.~\cite{Wi}) is the group  
denoted $\Hom(\AdsB,C)\rtimes\SpABd{R}$, with underlying set 
$\Hom(\AdsB,C)\times\SpABd{R}$ and product defined by:
$$
(\eta_1,\bar \phi_1)\cdot(\eta_2,\bar \phi_2)=(\eta_1 + \eta_2\circ \bar \phi_1\inv, \bar \phi_1\circ\bar \phi_2)
$$
It fits into a canonical exact sequence with the obvious homomorphisms
\begin{equation}\label{eq:ses for semidirect prod}
	1\rightarrow \Hom(\AdsB,C)
	\longrightarrow \Hom(\AdsB,C)\rtimes\SpABd{R}
	\longrightarrow
	\SpABd{R}\rightarrow 1
\end{equation}
%
for which the non-normal subgroup $(0,\SpABd{R})$ provides a natural splitting. 
\begin{thm}\label{thm:about Theta}
	The map 
	$$
	\displaymapdef{\Theta_\cD}{\Aut^0_R(\cH)}
	{\Hom(\AdsB,C)\rtimes\SpABd{R}}
	{\phi}
	{(\eta_{\phi},
	\bar{\phi}\inv)}
	$$
	is a group {\em anti}-isomorphism. 
\end{thm}
\bPf\ $\Theta_\cD$ is clearly a bijection, since $\t \Theta_\cD$ is. The anti-homomorphic property (\textit{viz} $\Theta_\cD(\phi_1\circ\phi_2)=\Theta_\cD(\phi_2)\cdot\Theta_\cD(\phi_1)$ for all $\phi_1,\phi_2\in\Aut^0_R(\cH)$) is proven by showing that 
$$
\eta_{\phi_1\circ\phi_2}=\eta_{\phi_2}+\eta_{\phi_1}\circ \bar{\phi}_2
$$
which follows from~(\ref{eq:phi(D(a))}). Details are left to the reader.\ePf
\rem\ This seems the most natural result, but if an isomorphism is required, one can simply compose $\Theta_\cD$ with inversion.\bigskip\\
Now consider diagram~(\ref{diag: 3x2 diagram of aut gps}) where the exact top row is by definition of $\Out^0_R(\cH)$, the bottom one is~(\ref{eq:ses for semidirect prod}) and we define the map $\Theta_1$  by 
$\Theta_1(\phi_h):=\eta_{\phi_h}$ for all $h\in\cH$.
\begin{figure}[ht]
	\begin{equation}\label{diag: 3x2 diagram of aut gps}
		{\xymatrix{
				1\ar[r]&{\rm Inn}(\cH)\ar[dd]^{\Theta_1}\ar[r]&\Aut^0_R(\cH)
				\ar[dd]^{\Theta_\cD}\ar[r]&
				\Out^0_R(\cH)
  				\ar@{-->}[dd]^{\t \Theta}\ar[r]&1\\
				&&&&&\\
				1\ar[r]&\Hom(\AdsB,C)   \ar[r]&\Hom(\AdsB,C)\rtimes\SpABd{R}\ar[r]&
				\SpABd{R}\ar[r]&1
		}}	
	\end{equation}
\end{figure}
Since 
$\bar{\phi}_{h}=\id_{\AdsB}$, the left-hand square commutes and so  $\Theta_1$ must be an injective anti-homomorphism (but see Remark~\ref{rem:first anti-iso is an iso and splitting}~\ref{part: Comments1} below). 
Moreover  
$$|\Hom(\AdsB,C)|=|\Hom(A,C)|\,|\Hom(B,C)|=|A|\,|B|=|\AdsB|=|\bar \cH|=|\Inn(\cH)|$$
so $\Theta_1$ is an anti-\textit{isomorphism}. Hence there is a unique anti-isomorphism   
$\t\Theta$ making the right square commute and since  
$
\t \Theta(\hat \phi)=\bar \phi\inv
$
we see that $\t \Theta$ is just the homomorphism $\bar \Theta$ of Proposition~\ref{prop:the hom bar Theta} composed with inversion.
We have proved the following:
 \begin{thm}\label{thm: 3x2 diagram of aut gps}
 Suppose condition ND-C and Hypothesis~\ref{hyp:r is odd} hold. Then diagram~(\ref{diag: 3x2 diagram of aut gps}) commutes and has exact rows. The vertical maps (described above) are anti-isomorphisms. Moreover, $\bar \Theta$ is an isomorphism.\ePf 
 \end{thm}
\rem\ {\bf (Comments on Theorem~\ref{thm: 3x2 diagram of aut gps})}
\label{rem:first anti-iso is an iso and splitting}
\begin{enumerate}
\item\label{part: Comments1}{\bf (The Map $\Theta_1$)}\ The groups $\Inn(\cH)\cong\cH/Z=\AdsB$ and $\Hom(\AdsB,C)$ are always \textit{abelian} and have natural $R$-module structures such that the action of $\SpABd{R}$ on the latter is $R$-linear. In particular, the map $\Theta_1$ is both an isomorphism and an anti-isomorphism. 
Note also that  equation~(\ref{eq:phi(D(a))}) for $\phi=\phi_h$ reads
$$
h\cD(\ba)h\inv=\cD(\ba)m(\eta_{\phi_h}(\ba))
\addwords{for all $\ba\in \AdsB$}
$$ 
Thus 
$
m(\eta_{\phi_h}(\ba))=[h,\cD(a)]=\delta_\cH(\bar{h},\ba)
$
and so  $\Theta_1(\phi_h)=\eta_{\phi_h}$ is simply the homomorphism $\delta(\bar{h},\cdot)$ in $\Hom(\AdsB,C)$. This definition of $\Theta_1$ clearly makes sense without Hypothesis~\ref{hyp:r is odd}. It is always an $R$-isomorphism because $\delta$ is non-degenerate and $R$-balanced, 
\item\label{part: Comments2} {\bf (The Subgroup $\Aut_R^0(\cH)_\cD$)}\ 
The natural splitting of the bottom row of~(\ref{diag: 3x2 diagram of aut gps}) induces one of the top row \textit{via} $\Theta_\cD$, namely that associated to  the non-normal subgroup of $\Aut^0_R(\cH)$ consisting of those automorphisms, $\phi$ such that $\eta_\phi$ is trivial, \ie\ those that `commute with $\cD$' in the sense that  $\phi(\cD(\ba))=\cD(\bar{\phi}(\ba))$ for all $\ba$. We therefore denote it $\Aut_R^0(\cH)_\cD$. It obviously intersects $\Inn(\cH)$ trivially and maps isomorphically onto $\Out_R^0(\cH)$.  Notice also that $\Aut_R^0(\cH)_\cD$ contains the subgroup $\im(\Delta)$ of diagonal $R$-automorphisms  (\cf\ Corollary~\ref{cor:Diagonal Automorphisms})). Indeed, it follows from the definitions of $\cD$ and  $\Delta$ that,
for any $\alpha\in\Aut_R(A)$ and  $a\in A, b\in B$, we have 
\begin{eqnarray}
	\Delta(\alpha)(\cD((a,b)))
	&=&h\left(\alpha(a),\beta(b),\frac12 \lambda(a,b)\right)\nonumber\\
	&=&h\left(\alpha(a),\beta(b),\frac12 \lambda(\alpha(a),\beta(b))\right)\nonumber\\
	&=&
	\cD(\alpha(a),\beta(b))\label{eq:D commutes w diag autoorphisms}
\end{eqnarray}
where the middle equality comes from the definition~(\ref{eq: relation between alpha and beta in diag aut}) of $\beta$.
\item\label{part: Comments3} 
The theorem leads to a description of $\Aut_R(\cH)$ itself. Indeed the sequence 
$$
\ses{1}{\Aut_R^0(\cH)}{\Aut_R(\cH)}{\Aut(C)=(\Zover{r})^\times}{1}
$$
is clearly exact on the left and in the middle. It is also exact on the right and splits: since $r$ kills $A$ and $B$ we may consider (for example) the splitting homomorphism from $(\Zover{r})^\times$ to $\Aut_R(\cH)$ sending $\bar i$ to the automorphism $\bar i \times \id_B \times \bar i$ for any $i\in\bbZ$ with $(i,r)=1$. 
\item\label{part: Comments4}
Finally, for versions of the theorem and the above comments `without subscript $R$'s' we can obviously just take $R=\bbZ$! 
\end{enumerate}
\noindent The study of the subgroups of  $\Out^0_R(\cH)$ now  reduces to that of the subgroups of $\SpABd{R}$. The latter can be made explicit when, for instance, $\AdsB$ is a free $R$-module, as in the examples of Section~\ref{sec:Examples}.
\subsection{Displacement Operators in the Schr\"odinger Representation}\label{subsec:Projectors}
Now suppose $p:C\rightarrow\bbC^\times$ is an injective homomorphism and define the homomorphism $\psi:Z\rightarrow\bbC^\times$ by $\psi\circ m = p$ so that $\sigma_p$ is a faithful SV-representation of $\cH$ on $\cM(A)$ with character $\chi_p$ and central character $\psi$ (by Proposition~\ref{prop:faithful SV reps for HeiStan under Hyp 5.1}).
We define the {\em displacement operators} (w.r.t.\ $p$) by  
$$
D(\ba)=D_p(\ba):=\sigma_p(\cD(\ba))\in{\rm U}(\cM(A))\addwords{for all $\ba\in A\oplus B$}
$$
Applying $\sigma_p$ to equation (\ref{eq:D is a hom mod Z}) shows in particular that the map $D:\ba\mapsto D(\ba)$ induces a homomorphism $A\oplus B\rightarrow {\rm PU}(\cM(A))$. (It is not hard to show that this cannot be lifted to a homomorphism of $A\oplus B$ into ${\rm U}(\cM(A))$.)
The adjoint $D(\ba)^\dagger$ w.r.t. $\langle\ ,\ \rangle$ is equal to $D(\ba)\inv=\sigma_p(\cD(\ba)\inv)=D(-\ba)$ (using~(\ref{eq:D is a hom mod Z}) for $n=-1$).
\begin{prop} With the above assumptions and notations  
	$$\Tr(D(\ba)^\dagger D(\ba'))=\eitherortw{s}{\mbox{if $\ba=\ba'$ in $A\oplus B$}}{0}{otherwise} $$
	and the subset ${\bD}:=\{D(\ba)\,:\,\ba\in A\oplus B\}$ of $   {\rm U}(\cM(A))$ is a $\bbC$-basis for $\End_\bbC(\cM(A))$.
\end{prop}
\bPf\ Note that $\Tr(D(\ba)^\dagger D(\ba'))$ equals $\chi_p(\cD(\ba)\inv\cD(\ba'))$ which is zero unless $\cD(\ba)\inv\cD(\ba')$ lies in $Z$ (since $\sigma_p$ is SV) in which case $\ba=\ba'$ so $\cD(\ba)\inv\cD(\ba')=\id_\cH$. The equation follows. An easy argument then shows that $\bD$ is linearly independent and since $|\bD|=s^2$, it is a basis of $\End_\bbC(\cM(A))$. (Alternatively, use Theorem~\ref{thm:equivalent condits for SV}~\ref{part4:equivalent condits for SV}.)\ePf 
\noindent 
\rem\ The proposition amounts to the statement that the normalised displacement operators $s^{-\frac12}D(\ba)$ for $\ba\in A\oplus B$ form an orthonormal basis of $\End_\bbC(\cM(A))$ w.r.t.\ the {\em Hilbert-Schmidt inner product}.\bigskip\\
\rem\label{rem:  angle-maps related to projectors} ({\bf Relation between the Overlap-Map $\cAng_{\cD,\ell}$ and `Physical' Overlaps} )\\
Let  $\ell\in \bbP\cM(A)$ be a line generated by a unit vector $v\in\cM(A)$, say, and let $\Pi_\ell\in\End_\bbC(\cM(A))$ be the orthogonal projector of $\cM(A)$ onto $\ell$. By the proposition, we have:
$$
\Pi_\ell=\sum_{\ba\in A\oplus B}\fl_{\ba}D(\ba)
\ \ \ \mbox{where}\ \ \  s\,\fl_{-\ba}=\Tr(D(\ba)\Pi_\ell)
$$
Note also that the matrix of $\Pi_\ell$ w.r.t.\ any orthonormal basis $\cB$ of $\cM(A)$ is just $\bv\bv^{T,\ast}$ where $\bv$ denotes the column vector of $v$ w.r.t.\ $\cB$.\bigskip\\
In the very special Base Case, the physics literature on SICs defines the \textit{overlaps} of $\ell$ or, equivalently, of $\Pi_\ell$, to be precisely the quantities $s\,\fl_\ba$ above. (\textit{Warning:} in that case, $\ba$ is naturally an element of $(a,b)\in(\Zover{d})^2$ yet -- for reasons made clearer in~\S\ref{subsubsection:base case odd d} -- the complex number corresponding to our $s\,\fl_{(a,b)}$ is conventionally labelled with the pair $(-a,b)$.)
To relate these to the \textit{overlap-maps} for Generalised Heisenberg Groups in our set-up, we can calculate $\Tr(D(\ba)\Pi_\ell)$ by extending $\{v\}$ to an orthonormal basis $\cV:=\{v_1,v_2,\ldots v_s\}$ of $\cM(A)$ with $v_1=v$. Indeed, 
the coefficient of $v_i$ in $D(\ba)\Pi_\ell v_i$ is just $\langle v,D(\ba)v\rangle$ if $i=1$ and otherwise $0$. So
$
\Tr(D(\ba)\Pi_\ell)=\langle v,D(\ba)v\rangle
$ 
Thus, the definition~(\ref{eq:defn of Overlap-Map}) of our overlap-map (still identifying $\bar\cH$ with $\AdsB$) gives 
\begin{equation}\label{eq:connection between overlap maps and l's}	
\cAng_{\cD,\ell}(\ba)=\langle v,\cD(\ba)\cdot_{\sigma_p}v\rangle=
\langle v,D(\ba)v\rangle=
s\,\fl_{-\ba}
\end{equation}
It follows more generally that 
$
\langle v,h(\ba,c)\cdot_{\sigma_p}v\rangle=\psi\left(c-\frac12\lambda(a,b)\right)s\,\fl_{-\ba}
$.
\section{GHGs of Arithmetic Type: an Introduction}\label{sec:Examples} 
\subsection{Definitions} 
In the general set-up, we take $k$ to be be any number field (of degree $n\geq 1$ over $\bbQ$, say) and $\cO$ to be any order in $k$. For the purposes of this introduction, however, we shall consider only the case  where \textit{all non-zero ideals of $\cO$ are invertible.} Equivalently, $\cO$ is maximal, \ie\  $\cO=\cO_k$, the ring of algebraic integers in $k$. We choose $I\subset k$ to be be a fractional $\cO$-ideal with inverse $I\inv$ and write $\hat{I}$ for the fractional ideal $\fD\inv I\inv$ where $\fD=\fD_{k/\bbQ}$, is the absolute different of $k$ (an integral ideal of $\cO$, see \cite[Ch.\ 3]{Serre Local fields}). 
%
Now suppose that $\ff$ is a non-zero, proper ideal of $\cO$ and define an integer $f>1$ by $\ff\cap\bbZ=f\bbZ$, so that $\ff|f\cO$. Then $I/\ff I$ and  $\widehat{\ff I}/\hat I=\ff\inv\hat I/\hat I$ are natural $\cO/\ff$-modules. By the standard theory of the different, the absolute trace map  $\Tr=\Tr_{k/\bbQ}:k\rightarrow \bbQ$ now gives rise to a well-defined, bilinear, $\cO/\ff$-balanced, non-degenerate pairing
$$
\displaymapdef{\lambda=\lambda_{I,\ff}}{I/\ff I\times \ff\inv  \hat I/\hat{I}}
{f\inv\bbZ/\bbZ}
{(\bar{a},\bar{b})}
{\Tr(ab)+\bbZ\addwords{for all $a\in I$ and $b\in\ff\inv\hat I$}}
$$
\rem\label{rem:two simplifications}\ 
The ideals $\ff$ and $\fD$ are relatively prime iff all prime ideals of $\cO$ dividing $\ff$ are unramified over $\bbZ$. In this situation $\fD+\ff=\cO$ so that $\ff\inv I\inv+\hat I=\ff\inv I\inv+\fD\inv I\inv=\ff\inv\fD\inv I\inv=\ff\inv\hat I$. It follows that the inclusion $\ff\inv I\inv\hookrightarrow \ff\inv \hat I$ induces an 
isomorphism from $\ff\inv I\inv/I\inv$ to $\ff\inv  \hat I/\hat{I}$, allowing us to replace the latter by the former in the definition of $\lambda_{I,\ff}$. 
\bigskip\\ 
To simplify notation we set 
$
A_{I,\ff}:= I/\ff I$ and $B_{I,\ff}:= \ff\inv \hat I/\hat I
$
and define 
$$ 
\cH[I,\ff]:=\Hei{A_{I,\ff}}{B_{I,\ff}}{f\inv\bbZ/\bbZ}{\lambda_{I,f}}
$$ 
It is a generalised Heisenberg group with $\Of$-structure, satisfying condition ND-C. We call it a GHG \textit{of arithmetc type}.  
In the notation of Subsection~\ref{subsec:ND-C}, the integer $s$ equals ${\rm N}\ff$ (the norm) and the integers $e$ and $r=|Z(\cH[I,\ff])|$ are both equal to $f$.\bigskip\\
\rem\label{rem:enlarging the centre} {\bf(Enlarging the Centre)} It is sometimes useful to enlarge the centre of $\cH[I,\ff]$ by setting 
$$
\cH[I,\ff,r]:=\Hei{A_{I,\ff}}{B_{I,\ff}}{r\inv\bbZ/\bbZ}{\lambda_{I,f}}
$$ 
for \textit{any} positive multiple $r$ of $f$. Thus $\cH[I,\ff,r]$ has centre isomorphic to $\Zover{r}$ and a subgroup of index $r/f$ isomorphic $\cH[I,\ff]$ with the  same quotient-by-centre. (\textit{Warning}: in earlier work, the notation `$\cH[I,\ff,r]$' was used to denote  
$\Hei{I/\ff I}{r\ff\inv\hat I/r\hat I}{\bbZ/r\bbZ}{\bar{\Tr}}$ which is, however, isomorphic to the current definition of $\cH[I,\ff,r]$ \textit{via} the diagonal isomorphism $\id \times r \times r$, in an obvious notation.)\bigskip\\
Consider the  injective homomorphism
$\be:\bbQ/\bbZ\rightarrow\bbC^\times$ sending $t+\bbZ$ to $\exp(2\pi it)$. Taking $p$ in Proposition~\ref{prop:faithful SV reps for HeiStan under Hyp 5.1} to be the restriction of $\be$ to $\Zunder{f}$ gives a faithful, irreducible, unitary  Schr\"odinger representation $\sigma_\be$ of $\cH[I,\ff]$ on $\cM(A_{I,\ff})$ (the vector space of complex functions on $I/\ff I$) defined by equation~(\ref{eq: defn of Schrod psi-action on functions}) with $\lambda=\lambda_{I,\ff}$. Equivalently, we have an injective homomorphism  
$$
\sigma_\be:\cH[I,\ff]\longrightarrow \U(\cM(I/\ff I))
$$
\subsection{The Automorphism Group of $\cH[I,\ff]$}
Abbreviating $\cH[I,\ff]$ to $\cH$, we identify $\cH/Z(\cH)$ with $A_{I,\ff}\oplus B_{I,\ff}$ as usual. 
For convenience, we shall write its elements as column vectors, so, by definition,  
\begin{equation}\label{eq:def delta now}
	\delta((\overline{a}, \overline{b})^T,(\overline{a'},\overline{b'})^T)=
	\overline{\Tr(ab'-a'b)}\in \Zunder{f}
	\addwords{for all $a,a'\in I$, $b,b'\in \ff\inv I$}
\end{equation}
By Proposition~\ref{prop:the hom bar Theta}, we have well-defined homomorphism 
$\bar\Theta:\Out^0_{\Of}(\cH)\rightarrow{\rm Sp}_{\Of}(A_{I,\ff}\oplus B_{I,\ff};\delta)$ sending $\hat\phi$ to $\bar \phi$ for all $\phi\in\Aut^0_\Of(\cH)$. In the case 
that $r$ is \textit{odd} Theorem~\ref{thm: 3x2 diagram of aut gps} says it is an \textit{isomorphism}. 
We now exhibit an explicit (but largely non-canonical)
 isomorphism from ${\rm Sp}_{\Of}(A_{I,\ff}\oplus B_{I,\ff},\delta)$ to 
 $\SL_2(\Of)$. 
 Similar analyses can be performed for ${\rm Sp}_R(A_{I,\ff}\oplus 
 B_{I,\ff},\delta)$, considering $A_{I,\ff}$ and $B_{I,\ff}$ as modules over other choices of ring $R$, \eg\  $R=\bbZ$. However, it is particularly simple for $R=\Of$ 
because  $A_{I,\ff}$ and $B_{I,\ff}$ are both \textit{free of rank 1} over this ring. 
In fact, the Weak Approximation Theorem allows us to choose $x\in I$ such that $x\cO+\ff I=I$. Thus we have surjective $\Of$-module homomorphism from $\Of$ onto $I/\ff I=A_{I,\ff}$ sending $a\in\Of$ to $a\bar x$. Since $|\Of|=N\ff=|I/\ff I|$, it must be an isomorphism, \ie\  $\{\bar x\}$ is an $\Of$-basis  for $A_{I,\ff}$. In the same way, we choose $y\in \ff\inv \hat I$ such that $y\cO+\hat I=\ff\inv \hat I$ and $\{\bar y\}$ is then an $\Of$-basis for $B_{I,\ff}$. There are no canonical choices of $\bar x$ and $\bar y$ in general although, in certain situations the two can be linked. 
Nevertheless, writing elements of $\AIdsB$ as column vectors, we have an ordered $\Of$-basis  
 $
%
((\bar x,0)^T,(0,\bar y)^T)
 $  
 of $\AIdsB$ and so an isomorphism of $\Of$-algebras
\begin{equation}\label{eq:def Xi}
\displaymapdef{\Xi=\Xi_{\bar x,\bar y}}
{\End_\Of(\AIdsB)}{{\rm M}_2(\Of)}{\bar \phi}{M_{\bar \phi}:=
\smallttmat{\bar u_{\bar \phi}}{\bar v_{\bar \phi}}{\bar w_{\bar \phi}}{\bar z_{\bar \phi}}
}
\end{equation}
where $u_{\bar\phi}, v_{\bar\phi}, w_{\bar\phi}$ and $z_{\bar\phi}$ are elements of $\cO$ whose images  $\barphi{u},\barphi{v},\barphi{w},\barphi{z}$ in $\Of$ are uniquely defined by  
$\bar \phi((\bar x,0)^T)=(\barphi{u}\bar x,\barphi{w}\bar y)^T$ and    
$\bar \phi((0,\bar y)^T)=(\barphi{v}\bar x,\barphi{z}\bar y)^T$. 
Restricting $\Xi$ to unit-groups gives an isomorphism 
$\Xi^\times=\Xi^\times_{\bar x,\bar y}$ from $\Aut_\Of(\AIdsB)$ to $\GL_2(\Of)$. 
\begin{thm}\label{thm: Sp iso to SL in general arithmetic case} For each pair of elements $\bar x$ and $\bar y$ as above, the isomorphism 
	$\Xi^\times$ restricts further to an isomorphism from 
	${\rm Sp}_{\Of}(A_{I,\ff}\oplus B_{I,\ff},\delta)$ to 
	$\SL_2(\Of)$.
\end{thm}
\bPf\ By definitions~(\ref{eq:defn SpABd}), ${\rm Sp}_{\Of}(A_{I,\ff}\oplus B_{I,\ff},\delta)$ consists of the elements $\bar \phi$ of $\Aut_\Of(\AIdsB)$ satisfying 
the equality 
\begin{equation}\label{eq:when phi is in Sp}
\delta(\bar\phi(\bx),\bar\phi(\bx'))=\delta(\bx,\bx')
\end{equation}
for all $\bx,\bx'\in\AIdsB$. Since both sides are $\bbZ$-bilinear, $\Of$-balanced and skew-symmetric as functions of $(\bx,\bx')$, this is equivalent to~(\ref{eq:when phi is in Sp}) holding for $(\bx,\bx')$ equal to $((\bar c\bar x,0)^T, (\bar x,0)^T)$,  
$((0,\bar c\bar y)^T, (0,\bar y)^T)$ and $((\bar c\bar x,0)^T,(0,\bar y)^T)$, for any $c\in \cO$. In the first two cases, one checks easily using~(\ref{eq:def delta now}) that equation~(\ref{eq:when phi is in Sp}) reads `$0=0$' for any $\bar\phi$. In the third case, it reads
$$ 
\delta((\overline{cu_{\bar \phi}x},
\overline{cw_{\bar \phi}y})^T,
(\overline{v_{\bar \phi}x},\overline{z_{\bar \phi}y})^T)=
\delta((\overline{cx},0)^T,(0,\overline{y})^T)
$$
where we are using the notation of equation~(\ref{eq:def Xi}) and following. Thus, using~(\ref{eq:def delta now}), we see that $\bar\phi$ 
lies in  ${\rm Sp}_{\Of}(A_{I,\ff}\oplus B_{I,\ff},\delta)$ iff
$$
\overline{\Tr((u_{\bar \phi}z_{\bar \phi}-v_{\bar \phi}w_{\bar \phi})cxy)}=\overline{\Tr(cxy)}
\addwords{in $\Zunder{f}$ for all $c\in\cO$}
$$ 
or, in other words, iff
$$
\Tr((d_{\bar \phi}-1)cxy)\in \bbZ
\addwords{for all $c\in\cO$}
$$
where $d_{\bar \phi}:=u_{\bar \phi}z_{\bar \phi}-v_{\bar \phi}w_{\bar \phi}\in\cO$ or, equivalently, iff 
 \begin{equation}\label{eq: crit for phi}
 (d_{\bar \phi}-1)xy\in \fD\inv
 \end{equation}
I claim that equation~(\ref{eq: crit for phi}) holds iff $d_{\bar \phi}-1$ lies in $\ff$. Since  $\det (M_{\bar\phi})=\overline{d_{\bar \phi}}$, this in turn is equivalent to $\det(\Xi^{\times}(\bar\phi))=\bar{1}$ in $\Of$ and the result will follow. But $d_{\bar \phi}-1\in\ff$ 
certainly implies~(\ref{eq: crit for phi}), since $xy\in \ff\inv I\hat I=\ff\inv \fD\inv$. For the converse implication, note that,  by choice of $x$ and $y$, we have 
$$
\ff\inv\fD\inv=I (\ff\inv \hat I) =xy\cO+x\hat I+y\ff I+\ff \fD\inv \subset xy\cO + \fD\inv
$$
Multiplying this containment through by $d_{\bar \phi}-1\in\cO$, we see that 
equation~(\ref{eq: crit for phi}) implies $(d_{\bar \phi}-1)\ff\inv \fD\inv \subset \fD\inv$ so $d_{\bar \phi}-1\in \ff$, proving the claim, and hence the theorem.\ePf 
\rem\ {\bf (Variation of $\bar x$, $\bar y$,  and Diagonal Matrices)} Changing the basis elements $\bar x$ and $\bar y$ will obviously compose $\Xi^\times_{\bar x, \bar y}$ with conjugation by a \textit{diagonal} transition matrix in $\GL_2(\Of)$. Since $A_{I,f}$ is free of rank-1, $\Aut_\Of(A_{I,f})$ is {\em canonically} identified with $(\Of)^\times$ and  for every $u$ in the latter, the \textit{diagonal automorphism} $\Delta(u)=u\times u\inv \times \id_{\Zunder{f}}$ of $\cH$ 
is mapped by the composite homomorphism 
$$
\Psi_{\bar x, \bar y}:\Aut^0_\Of(\cH)
\longrightarrow
\Out^0_\Of(\cH)
\stackrel{\bar \Theta}{\longrightarrow}
{\rm Sp}_{\Of}(A_{I,\ff}\oplus B_{I,\ff},\delta)
\stackrel{\Xi_{\bar x, \bar y}^\times}{\longrightarrow}
\SL_2(\Of)
$$
to the \textit{diagonal matrix} ${\rm diag}(u,u\inv)$, irrespective of the choice of $\bar x$ and $\bar y$.
\subsection{Perspectives}\label{subsec:perspectives} In a future paper, we shall develop the following aspects of GHGs of arithmetic type, $\cH[I,\ff]$, their Schr\"odinger and Weil representations, bouquets and $p$-adic limits:
\begin{enumerate}
	\item the natural relation between $\cH[I,\ff]$ and 
	$\cH[I',\ff]$ and their respective bouquets when $I$ and $I'$ lie in the same ideal class.
	\item given an ideal  $\fg$ dividing $\ff$, two related processes of `descent' that allow us to represent both $\cH[I,\fg]$ and $\cH[I\fg,\fg\inv\ff]$ as different \textit{subquotients} of $\cH[I,\ff]$. We shall also study the corresponding effects on Schr\"odinger representations and bouquets, 
	\item  explicit versions of the (projective) Weil representation of $\Aut^0_\Of(\cH[I,\ff])$ and also of an affine lifting over the subgroup $\Aut^0_\Of(\cH[I,\ff])_\cD\cong \SL_2(\cO/\ff)$ in certain cases. 
	\item\label{part:rem further directions1} (for an odd prime $p$), two notions of $p$-adic limit for the groups $\cH[I,p^n\cO]$ as $n\rightarrow \infty$. Schr\"odinger representations of these limits are naturally on spaces of \textit{$\bbC_p$-valued} measures on $I\otimes\bbZ_p$ (on which the unit group $E_k$ acts naturally). 
	\item for one of these limits -- a compact $p$-adic Heisenberg group with centre isomorphic to $\bbZ_p$ -- we shall consider  an extension of the natural Schr\"odinger action defined via $p$-power roots of unity, to actions involving the continuous $p$-adic central character   $\psi_\kappa:\bbZ_p\rightarrow \bbC_p^\times$ , \textit{not necessarily of finite order}, which sends $t\in\bbZ_p$ to $\kappa^t$ for $\kappa$ in the `open' unit disc about $1$ in $\bbC_p$. Crucially, in the case $\kappa$ is a primitive $p^n$th root of unity, this descends back to the Schr\"odinger representation of $\cH[I, p^n\cO]$.   
	\item for general $\kappa$, we shall prove a result that substitutes for a sort of $p$-adic `Weil Representation' for the above compact limit, by means of certain `theta-like' $p$-adic double integrals associated to binary quadratic forms. 
\end{enumerate}
We hope that in this way, a combination of algebraic and $p$-adic analytic techniques may eventually lead to class-field theoretic applications and illuminate 
some of the remaining mysteries of SIC phenomenology in certain cases. Indeed, we  conclude this paper by explaining the connection between SICs and GHGs of arithmetic type in the simplest possible case.  
\subsection{The Base Case and Heisenberg SICs}\label{subsec:the Base Case} 
\subsubsection{The Base Case for $d$ Odd}\label{subsubsection:base case odd d}
Let $k=\bbQ$, $I=\bbZ$ and $\ff=d\bbZ$ for an \textit{odd} integer $d>2$. Thus $f=d$ and 
$$
\cH[\bbZ,d\bbZ]=\cH(\Zover{d},\Zunder{d},\Zunder{d},\times)
$$  
is a GHG of arithmetic type, of order $d^3$ and exponent $d$, nilpotent of class $2$. The bilinear form is just the obvious multiplication
$(a+d\bbZ, b+\bbZ)\mapsto ab+\bbZ$ for any $a\in\bbZ$, $b\in d\inv\bbZ$. There is an obvious diagonal isomorphism $\id\times d\times d$ mapping $\cH[\bbZ,d\bbZ]$ to 
$\cH(\Zover{d},\Zover{d},\Zover{d},\times)$. For notational reasons, we shall work mainly with the latter group, denoted $\cH_d$. Its group-law  is given explicitly by  
$$
h(a,b,c)h(a',b',c')=h(a+ a',b + b',c +c '+ab')
\addwords{for all $a,a',b,b'c,c'\in\Zover{d}$}
$$
(It follows that $\cH_d$ may also be represented as the group of $3\times 3$ upper-triangular, unipotent matrices with coefficients in $\Zover{d}$, \textit{cf}~\cite{Wi: Heisenberg Gp}.) We now unwind various definitions and results of the preceding sections and relate them to the published work on Heisenberg SICs in $\bbC^d$ mentioned the Introduction.

First, we have  $Z_d:=Z(\cH_d)=m(\Zover{d})=h(0,0,\Zover{d})\cong{\Zover{d}}$ and we may identify $\bar{\cH}_d=\cH_d/Z_d$ with $(\Zover{d})^2$. 
We set $h_1:=h(-\bar 1,0,0)$ and $h_2:=h(0,\bar 1,0)$, so that $[h_2,h_1]=h(0,0,\bar 1)$ generates $Z_d$. It follows easily that $h_1$ and $h_2$ generate $\cH_d$. 

Next, the isomorphism from $\cH[\bbZ,d\bbZ]$ to $\cH_d$ makes the Schr\"odinger representation $\sigma_\be$ defined for the former correspond to the faithful, $d$-dimensional SV representation
$
\sigma_d: \cH_d \longrightarrow 
\U(\cM(\Zover{d})) 
$
given by the explicit action 
\begin{equation}\label{eq:Schrod rep in Base Case1}
	h(a,b,c)\cdot_{\sigma_d} f(x)=\zeta_d^{bx+c}f(x+a)
	\addwords{for all $a,b,c,x$ in $\Zover{d}$ and $f\in\cM(\Zover{d})$}
\end{equation}
where $\zeta_d:=\exp(2\pi i/d)$. Thus, the 
central character $\psi_d$ sends $h(0,0,\bar c)$ to $\zeta_d^c$. 
Referring everything to the ordered basis $\cE_d:=(e_{\bar{0}},\ldots,e_{\overline{d-1}})$ of $\cM(\Zover{d})$ gives a isomorphism from $\cM(\Zover{d})$ to 
$\bbC^d$. Then $\langle\ ,\ \rangle$ becomes the standard PDHF on $\bbC^d$, and $\sigma_d$ becomes an injective homomorphism  
$\Sigma_d$ from $\cH_d$ to $\U(\bbC^d)$ acting naturally on $\bbC^d$. Using~(\ref{eq:Schrod rep in Base Case1}), we find explicitly:  
\begin{equation}\label{eq:defns of X and Z}
	\Sigma_d(h_1)=X:=
	\begin{pmatrix}
		0&0&\cdots&0&1\\	
		1&0&\cdots&0&0\\
		0&1&\cdots&0&0\\	
		\vdots&\vdots&\ddots&\vdots&\vdots\\	
		0&0&\cdots&1&0
	\end{pmatrix}
	\ \ \ \mbox{and}\ \ \ 
	\Sigma_d(h_2)=Z:=\begin{pmatrix}
		1&0&0&\cdots&0\\	
		0&\zeta_d&0&\cdots&0\\	
		0&0&\zeta_d^2&\cdots&0\\	
		\vdots&\vdots&\vdots&\ddots&\vdots\\	
		0&0&0&\cdots&\zeta_d^{d-1}
	\end{pmatrix}
\end{equation}
Thus $\Sigma_d$ induces an isomorphism from $\cH_d$  onto the 
subgroup of $U_d(\bbC)$ generated by the matrices $X$ and $Z$.
Since $[Z,X]=\Sigma_d(h(0,\bar 0, \bar 1))$ is $\zeta_d$ times the identity matrix, this subgroup is exactly 
the \textit{Weyl-Heisenberg Group} as defined for odd $d$ in~\cite{Kopp IMRN}, for example. It is variously denoted as ${\rm H}(d)$ in \textit{idem}, as ${\rm GP}(d)$ in~\cite{DMA: SICPOVMs and the ECG} (where it is called the \textit{Generalised Pauli Group}) and as ${\rm WH}(d)$ elsewhere, in the Introduction and here.

Now for the bouquets. Let $\ell\in \bbP \cM(\Zover{d})$ be a complex line in $\cM(\Zover{d})$ and let $\cY$ be the $\cH_d$-\textit{bouquet} which is its orbit under the action of $\bar \cH_d$ on $\bbP \cM(\Zover{d})$ induced by $\sigma_d$. 
We shall assume that  $\cY$ is free, \ie\ $|\cY|=d^2$.
Identifying $\bbP \cM(\Zover{d})$ with $\bbP \bbC^d$ using the basis $\cE_d$ as above, $\cY$ becomes the orbit under ${\rm WH}(d)$ of the line $l_0\in \bbP \bbC^d$ corresponding to $\ell$. Clearly $\cY$ is \textit{equiangular} iff this orbit is a (Heisenberg) SIC.

More generally, let $\cR: (\Zover{d})^2=\cH_d/Z_d\rightarrow \cH_d$ be any right-inverse to the quotient homomorphism and consider the overlap-map 
$
\cAng_{\cR,\ell}:(\Zover{d})^2\rightarrow\bbC
$
sending $(a,b)^T$ to $\langle v, \cR((a,b)^T)\cdot_{\sigma} v\rangle$ for any unit vector $v$ generating $\ell$ (see~(\ref{eq:defn of Overlap-Map}) and~(\ref{eq:connection between overlap maps and l's})). It descends to a well-defined map $\cAng_{\cY}:(\Zover{d})^2\rightarrow\mu_d\backslash\bbC$ which is independent of $\cR$, as in~(\ref{eq:overlap map of bouquet}), and depends  on $\ell$ only through $\cY$. Composing with the absolute value $|\cdot |$ gives the angle-map $\fa_{\cY}:(\Zover{d})^2\rightarrow[0,1]$
of the bouquet $\cY$. 
%
It follows from Theorem~\ref{thm:equiangular existence and clinometric} that $\cY$ will be equiangular iff $\fa_\cY(\gamma)^2=(d+1)\inv$ for all for all non-zero $\gamma\in(\Zover{d})^2$. 

Before discussing regular bouquets for $\cH_d$, we look more closely at its automorphism group and displacement operators which are well-defined since $d$ is odd. Indeed, writing elements of $\cH_d/Z_d=(\Zover{d})^2$ as column vectors $\ba=(a,b)^T$, equation~(\ref{def:displacement operator}) gives  
$\cD((a,b)^T)=h(a,b,\half ab)$ so that  $\Sigma_d(\cD((a,b)^T))$
is the element $\tau^{-ab}
X^{-a}Z^b
$
of ${\rm WH}(d)$ where $\tau=\zeta_d^{(d+1)/2}$. (This is denoted 
$D_{-a,b}$ in~\cite{Kopp IMRN}.) 
There is a 
{\em canonical} identification of $\Aut_{\Zover{d}}(\cH_d/Z_d)$ with $\GL_2(\Zover{d})$, sending $\alpha$ to the matrix $N_{\alpha}$, say. Since  $\delta(\ba,\ba')=ab'-a'b=\det(\ba|\ba')\in\Zover{d}$ , it is clear that the subgroup  $\Sp=\Sp_{\Zover{d}}(\Zover{d}\oplus \Zover{d})$ identifies canonically with the subgroup $\SL_2(\Zover{d})$. Thus  Theorem~\ref{thm:about Theta} 
establishes an\textit{ anti-}isomorphism 
\begin{equation}\label{eq:iso of Aut w  Z/d squared semidp SL Z/d in Base case}
	\displaymapdef
	{\Theta_d}
	{\Aut^0(\cH_d)}
	{\Hom((\Zover{d})^2,\Zover{d})\rtimes \SL_2(\Zover{d})}
	{\phi}
	{(\eta_\phi, N_{\bar \phi}\inv)}
\end{equation}
where $\eta_\phi$ is defined by equation~(\ref{eq:phi(D(a))}), that is:  $h(0,0,\eta_\phi((a,b)^T))=\phi(\cD(a,b)^T)\cD(N_{\bar \phi}(a,b)^T)\inv$. (In the semidirect product, an element $N$ of $\SL_2(\Zover{d})$ acts on the left on $\eta\in \Hom((\Zover{d})^2,\Zover{d})$ by precomposition with the left multiplication of column vectors by $N\inv$.) And $\Theta_d$ descends to an anti-isomorphism from $\Out^0(\cH_d)=\Out_{\Zover{d}}^0(\cH_d)$ to $\SL_2(\Zover{d})$ sending $\hat \phi$ to $N_{\bar \phi}\inv$.\bigskip\\
\rem\label{rem:TheCliffordGroupandTheExtCG}\ {\bf (The Clifford Group)}
This is is defined in \eg~\cite{DMA: SICPOVMs and the ECG} to be the normaliser of ${\rm H}(d)={\rm WH}(d)$  in $U_d(\bbC)$. Thus, by Remark~\ref{rem: re the Clifford}, it   identifies with the pre-image in $\U(\cM(\Zover{d}))$ of the isomorphic image of $\Aut^0(\cH_d)$ in ${\rm PU}(\cM(\Zover{d}))$ under the Weil representation  $\cT_{\sigma}$. Equation~(\ref{eq:iso of Aut w  Z/d squared semidp SL Z/d in Base case}) therefore essentially delivers the isomorphism $f$ in Theorem~1 of~\cite{DMA: SICPOVMs and the ECG} in the case of odd $d$. 
In fact, Appleby gives considerably more, namely an \textit{explicit characterisation} of the Weil representation. See, for example Lemma~4 of~\cite{DMA: SICPOVMs and the ECG} for the Base Case, and the detailed account in~\cite[Ch.\ 14]{Waldron's book}. 
The so-called \textit{Extended Clifford Group} of~\cite{DMA: SICPOVMs and the ECG} identifies in a similar way with the pre-image of $\t \cT_\sigma(\widetilde{\Aut^0}(\cH_d))$ in $\GL_\bbR(\cM(\Zover{d}))$, referred again to the basis $\cE_d$.  (See Remark~\ref{rem:R-lin extn of Weil}. Indeed the element $\hat J$ of the Extended Clifford Group corresponds to our $\bc$.)\bigskip\\
We now consider \textit{$\cA$-regular} bouquets for $\cA=\Out^0(\cH_d)$. 
It is well known (and not hard to see) that the orbits of $\SL_2(\Zover{d})$ acting on $(\Zover{d})^2$ are the just the disjoint subsets 
\begin{equation}\label{eq:orbits in Zoverdsquared}
O_j:=\{(a,b)^T\in (\Zover{d})^2\,:\,\mbox{\textit{additive} order of $(a,b)^T$ is $j$}\}\subset (\Zover{d})^2
\end{equation}
as $j$ runs through the set $S_d$, say, of positive divisors of $d$. Let $W_d$ denote the subspace of $\cM((\Zover{d})^2)$ consisting of (complex-valued) functions which are constant on each $O_j$. Then a free bouquet $\cY$ is $\Out^0(\cH_d)$-regular iff $\fa_{\cY}^2$ lies in $W_d$. To analyse the clinometric relation~(\ref{eq:aysquared is evec}) that $\fa_{\cY}^2$ must also satisfy, we write $\bw_j$ for the characteristic function of the subgroup $j(\Zover{d})^2<(\Zover{d})^2$ for each $j\in S_d$. One checks that $\{\bw_j:j\in S_d\}$ is a basis of $W_d$ and that
\begin{equation}\label{eq:Upsilon on bwj}
\Upsilon_p(\bw_j)=(d/j)^2\bw_{d/j}
\end{equation}
(To prove~(\ref{eq:Upsilon on bwj}), observe that, for fixed $\gamma'\in(\Zover{d})^2$, the function  $\gamma\mapsto p(\delta_{\cH_d}(\gamma,\gamma'))$ is a character on $(\Zover{d})^2$. By the bilinearity and non-degeneracy of $\delta_{\cH_d}$, and the injectivity of $p$, it is trivial on $j(\Zover{d})^2$ iff $\gamma'\in (d/j)(\Zover{d})^2$, \textit{etc.}, details LTR.) It follows without difficulty that $\Upsilon_p$ acts on $W_d$ and that its eigenspace  in $W_d$ for the eigenvalue $|\bar\cH_d|/d=d$ has the basis 
$$
\cB_d:=\left\{
\bu_j:=\frac{j^2}{j^2+d}\bw_j+\frac{d}{j^2+d}\bw_{d/j}\,:\ j\in \t S_d\,
\right\}
$$
where $\t S_d$ denotes the set $S_d\cap [1,\sqrt{d}]$, of cardinality  $\lceil|S_d|/2\rceil$. Thus the clinometric relation~(\ref{eq:aysquared is evec}) forces $\fa_\cY^2$ to be a (real) linear combination of $\cB_d$. Note also that $\bu_j(\gamma)$ lies in $[0,1]$ for all $j\in\t S_d$ and all $\gamma\in(\Zover{d})^2$ but equals $1$ for all $j$ iff $\gamma=(\bar 0,\bar 0)^T$.

If $d$ is prime then, $\fa_\cY^2$ must be  equal to $\bu_1$ because this is then the sole element of $\cB_d$, and so $\cY$ must correspond to a SIC. 
Otherwise $\cB_d$ has at least $2$ elements and so the interior  of its (real) convex hull in $W_d$, denoted $\langle\cB_d\rangle^0$, is \textit{uncountable} but any $f$ in $\langle\cB_d\rangle^0$ is presumably a potential candidate for $\fa_\cY^2$. Indeed, such an $f$ would certainly satisfy both the clinometric relation and the conditions imposed by Proposition~\ref{prop: properties of ov and ang}~\ref{part:properties of ov and ang0}, while causing $\cY$ to be $\Out^0(\cH_d)$-regular but not equiangular.
%

This situation seems incompatible with the regularity condition forcing the overlap map to take significant algebraic values, as the equiangularity condition apparently does for SICs (see the Introduction). It is possible, of course, that not all functions satisfying the clinometric relation are attainable as angle-maps of free orbits of $\cH_d$. Alternatively, some additional condition(s) on regular bouquets may be required for significant algebraicity to emerge as the desired generalisation of SIC phenomenology. For instance, an appropriate rational value of $\fa_{\cY}^2$ on each $O_j$ might be imposed (analogous to the values $1$ and $(d+1)\inv$ which are automatic for equiangular bouquets, see above). Natural candidates for such values are available and their 
numerical investigation - as well as that of regular bouquets more generally -- suggests itself as a sensible next step. Unfortunately, as for SICs themselves, the required computations are hugely intensive in all but the simplest examples. So far, they have provided only the meagrest of evidence.

\subsubsection{The Base Case for $d$ Even}
For even $d$ we just sketch the connection with 
the Weyl-Heisenberg groups of the SIC literature. 
The necessary changes to the analysis of automorphism groups and bouquets are left to the reader. The latter are mostly obvious generalisations  of the odd case and for the former, the relevant information on the corresponding Clifford group appears in~\cite{DMA: SICPOVMs and the ECG} alongside the odd case.

We again take $k=\bbQ$, $I=\hat I=\bbZ$ and $\ff=d\bbZ$ in the General Case, but now we enlarge the centre by a factor of $2$, taking $r=2d$ (so $4|r$) in 
Remark~\ref{rem:enlarging the centre}. Thus for $d$ \textit{even} we first consider the GHG
$
\cH[\bbZ,d\bbZ,2d]=\Hei{\Zover{d}}{\Zunder{d}}{\Zunder{(2d)}}{\times}
$
where the bilinear form is again multiplication, with image $\Zunder{d}$, now  considered as a subgroup of $\Zunder{(2d)}$. This time, we apply the diagonal homomorphism $\id\times 2d\times 2d$ which maps $\cH[\bbZ,d\bbZ, 2d]$ isomorphically onto 
$\cH(\Zover{d},2\bbZ/2d\bbZ,\Zover{2d},\times)$ and define $\cH_d$ to be the latter when $d$ is even. In particular $\cH_d$ has order $2d^3$ and centre $Z_d$ isomorphic to $\Zover{2d}$ with $\cH_d/Z_d\cong(\Zover{d})^2$ still. The elements $h_1=h(-\bar 1, 0,0)$, 
$h_2=h(0,\bar 2, 0)$, \textit{now supplemented by $h_3:=h(0,0,\bar 1)$} suffice to  generate $\cH_d$. 

Restricting $\be$ to $\Zunder{(2d)}$ gives rise to a faithful, SV, unitary,  Schr\"odinger representation $\sigma_\be$ of $\cH[\bbZ,d\bbZ,2d]$ on $\cM(\Zover{d})$. This corresponds to the representation
$\sigma_d $, say, of 
the isomorphic group $\cH_d$
given explicitly by the action 
\begin{equation}\label{eq:Schrod rep in Base Case2}
	h(a,b,c)\cdot_{\sigma_d} f(x)=\zeta_{2d}^{bx+c}f(x+a)
	\addwords{for all $a\in \Zover{d}$, $b\in 2\bbZ/2d\bbZ$, $c\in\bbZ/2d\bbZ$ and $f\in\cM(\Zover{d})$}
\end{equation}
where $\zeta_{2d}:=\exp(\pi i/d)$. Referring to the same basis $\cE_d$ of $\cM(\Zover{d})$ as before, we now get an injective homomorphism $\Sigma_d$ from  $\cH_d$ into $\U(\bbC^d)$ whose image is  
generated by the matrices $X=\Sigma_d(h_1)$ and 
$Z=\Sigma_d(h_2)$ (defined exactly as in~(\ref{eq:defns of X and Z})) \textit{together with} $\Sigma_d(h_3)$, which is $\zeta_{2d}$ times the identity matrix. Once again, this subgroup is exactly 
the \textit{Weyl-Heisenberg Group} ${\rm WH}(d)={\rm H}(d)$ as defined for even $d$ (see~\cite{Kopp IMRN}, for example).

\end{document}